\documentstyle[11pt,mst-stylefile,harvard,Bbb11a]{article}

\catcode`@=12
\newcommand{\bbZ}{{\Bbb Z}}
\newcommand{\bbR}{{\Bbb R}}

\newcommand{\bbC}{{\Bbb C}}
\newcommand{\bbL}{{\Bbb L}}

\newcommand{\bbQ}{{\Bbb Q}}

\renewcommand{\cite}{\citeyear}

\begin{document}

\title{Integral representations of periodic \\ and cyclic fractional stable motions
\thanks{ This research was partially supported by the NSF grant DMS-0102410 at Boston University.}
 }

\author{
 Vladas Pipiras
\\  University of North Carolina at Chapel Hill
\and
 Murad S.Taqqu
\\  Boston University
}

\bibliographystyle{agsm}

\maketitle


\begin{abstract}
Stable non-Gaussian self-similar mixed moving averages can be decomposed into several
components. Two of these are the periodic and cyclic fractional stable motions which are
the subject of this study. We focus on the structure of their integral representations
and show that the periodic fractional stable motions have, in fact, a canonical
representation. We study several examples and discuss questions of uniqueness, namely how
to determine whether two given integral representations of periodic or cyclic fractional
stable motions give rise to the same process.
\end{abstract}


\section{Introduction}
\label{s:preliminar}

{\it Periodic and cyclic fractional stable motions} ({\it PFSMs} and {\it CFSMs}, in
short) were introduced by Pipiras and Taqqu \cite{pipiras:taqqu:2003cy-fr} in the context
of a decomposition of symmetric $\alpha$-stable ($S \alpha S$, in short),
$\alpha\in(0,2)$, self-similar processes $X_\alpha(t)$, $t\in\bbR$, with stationary
increments having a mixed moving average representation
\begin{equation}\label{e:mma}
\{X_\alpha(t)\}_{t\in\bbR} \stackrel{d}{=}\left\{ \int_X\int_\bbR \Big(G(x,t+u) -
G(x,u)\Big) M_\alpha(dx,du)\right\}_{t\in\bbR},
\end{equation}
where $\stackrel{d}{=}$ stands for the equality of the finite-dimensional distributions.
Here, $M_\alpha$ is a symmetric $\alpha$-stable random measure with control measure
$\mu(dx)du$ (see Samorodnitsky and Taqqu \cite{samorodnitsky:taqqu:1994book}), $(X,{\cal
X},\mu)$ is a standard Lebesgue space (defined, for example in Appendix A of Pipiras and
Taqqu \cite{pipiras:taqqu:2003cy}) and
\begin{equation}\label{e:G_t-in-Lalpha}
\{G_t\}_{t\in\bbR}\in L^\alpha(X\times\bbR,\mu(dx)du),
\end{equation}
where
\begin{equation}\label{e:G_t}
    G_t(x,u) = G(x,t+u) - G(x,u),\quad x\in X,u\in\bbR.
\end{equation}

A process $X_\alpha$ is said to be self-similar with index of self-similarity $H>0$ if
for any $a>0$, $X_\alpha(at)$ has the same finite-dimensional distributions as
$X_\alpha(t)$. We will be interested in self-similar processes which have stationary
increments as well ($H$-sssi processes, in short). Note that a process $X_\alpha(t)$ of
the form (\ref{e:mma}) always has stationary increments. In order for the process
(\ref{e:mma}) to be self-similar, it is necessary to impose additional conditions on the
function $G$. There are various ways of doing this and different choices may give rise to
different types of processes. In order to understand the nature of the resulting
processes, one can associate the process $X_\alpha$ in (\ref{e:mma}) (or its kernel $G$)
to flows. Flows are deterministic maps $\psi_c$, $c>0$, satisfying $\psi_{c_1c_2} =
\psi_{c_1}\circ \psi_{c_2}$, $c_1,c_2>0$. One can then use the characteristics of the
flows to classify the corresponding $H$-sssi processes as well as to decompose a given
$H$-sssi process into sub-processes that belong to disjoint classes. Periodic flows are
such that each point of the space comes back to its initial position in a finite period
of time. Cyclic flows are periodic flows such that the shortest return time is positive
(nonzero).

In this work, we examine PFSMs and CFSMs in greater depth. We show that PFSMs can be
defined as those self-similar mixed moving averages having the representation
(\ref{e:mma}) with
\begin{equation}\label{e:cpfsm-X}
    X = Z\times [0,q(\cdot)),\quad \mu(dx) = \sigma(dz)dv
\end{equation}
where $x = (z,v)$, and
$$
 G(z,v,u) = b_1(z)^{[v + \ln |u|]_{q(z)}} \Big( F_1(z,\{v+\ln|u|\}_{q(z)})\,
u_+^\kappa + F_2(z,\{v+\ln|u|\}_{q(z)})\, u_-^\kappa\Big)
$$
\begin{equation}\label{e:cpfsm-G}
+\ 1_{\{b_1(z) = 1\}} 1_{\{\kappa = 0\}} F_3(z) \ln|u|,
\end{equation}
where $(Z,{\cal Z},\sigma)$ is a standard Lebesgue space, $b_1(z)\in \{-1,1\}$, $q(z)>0$
a.e.\ $\sigma(dz)$, $F_1,F_2:Z\times [0,q(\cdot)) \mapsto \bbR$, $F_3:Z\mapsto \bbR$ are
some functions. We used here the convenient notation
\begin{equation}\label{e:kappa}
    \kappa = H - \frac{1}{\alpha},
\end{equation}
where $H>0$ is the self-similarity parameter of the process $X_\alpha$. Thus, in
particular, $\kappa > -1/{\alpha}$. We also let
\begin{equation}
[x]_{a} = \max\{n\in\bbZ : na \le x \}, \quad \{x\}_{a} = x - a [x]_{a},\quad x\in
\bbR,a>0, \label{e:int-frac-parts}
\end{equation}
and suppose by convention that
\begin{equation}\label{e:pm-kappa=0}
    u_+^\kappa = 1_{(0,\infty)}(u),\quad u_-^\kappa = 1_{(-\infty,0]}(u),
\end{equation}
when $\kappa =0$. Since the representation (\ref{e:mma}) with $X$ defined as in
(\ref{e:cpfsm-X}) and $G$ defined as in (\ref{e:cpfsm-G}) characterizes a PFSM, it is
called {\it a canonical representation} for PFSMs. We will provide other canonical
representations for PFSMs as well. We are not aware of canonical representations for
CFSMs. CFSMs, however, do admit the representation (\ref{e:mma}) with (\ref{e:cpfsm-X})
and (\ref{e:cpfsm-G}) because they are PFSMs.

By using the representation (\ref{e:mma}) with (\ref{e:cpfsm-X}) and (\ref{e:cpfsm-G}),
we will generate and study various PFSMs and CFSMs. For these processes to be
well-defined, one must choose functions $F_1,F_2,F_3$ in (\ref{e:cpfsm-X}) so that the
functions $\{G_t\}_{t\in\bbR}$ in (\ref{e:G_t}) belong to the space $L^\alpha(Z\times
[0,q(\cdot))\times \bbR, \sigma(dz)dvdu)$. This is in general quite difficult (see
Section \ref{s:examples}). We will also address the uniqueness problem of PFSMs and
CFSMs, namely, how to determine whether two given PFSMs or CFSMs are different, that is,
when their finite dimensional distributions are not the same up to a constant.

The paper is organized as follows. In Section \ref{s:pfsm-cfsm}, we briefly recall the
definitions of PFSMs and CFSMs, and related notions. In Section \ref{s:canonical-pfsm},
we establish several canonical representations for PFSMs. In Section
\ref{s:represen-cyclic}, we discuss the representation problem for CFSMs. Examples of
PFSMs and CFSMs are given in Section \ref{s:examples}. Uniqueness questions are addressed
in Section \ref{s:unique}. In Section \ref{s:functionals}, we study some functionals
related to cyclic flows. Section \ref{s:proofs-main} contains the proofs of some results
of Section \ref{s:canonical-pfsm}.


\section{Periodic and cyclic fractional stable motions}
\label{s:pfsm-cfsm}

Periodic and cyclic fractional stable motions (PFSMs and CFSMs, in short) can be defined
in two equivalent ways (see Pipiras and Taqqu \cite{pipiras:taqqu:2003cy-fr}). The first
definition uses the kernel function $G$ in the representation (\ref{e:mma}) of stable
self-similar mixed moving average. Consider the sets
\begin{eqnarray}
C_P & = & \Big\{ x\in X: \exists\  c = c(x)\neq 1: G(x,cu) = b\, G(x,u + a)
+ d\ \ \mbox{a.e.}\ du\nonumber \\
& &\hspace{.7in} \mbox{for some}\ a = a(c,x), b=b(c,x)\neq 0,d=d(c,x)\in\bbR \Big\},
\label{e:pfsm-set}
\end{eqnarray}
and
\begin{equation}\label{e:cfsm-set}
    C_L = C_P\setminus C_F,
\end{equation}
where
\begin{eqnarray}
C_F & = & \Big\{ x\in X: \exists\ c_n = c_n(x)\to 1\ (c_n\neq 1): G(x,c_nu)
= b_n\, G(x,u + a_n) + d_n\ \mbox{a.e.}\ du\nonumber \\
& &\hspace{.4in} \mbox{for some}\ a_n = a_n(c_n,x), b_n=b_n(c_n,x)\neq
0,d_n=d_n(c_n,x)\in\bbR \Big\}. \label{e:mlfsm-set}
\end{eqnarray}
The sets $C_P$, $C_L$ and $C_F$ are called, respectively, the PFSM set, the CFSM set and
the mixed LFSM set.

\begin{definition}\label{d:pfsm-cfsm-kernel}
A $S \alpha S$, $\alpha \in (1,2)$, self-similar mixed moving average $X_\alpha$ having a
representation (\ref{e:mma}) is called
\begin{eqnarray*}
  \mbox{PFSM} &\mbox{if}& X=C_P, \\
  \mbox{CFSM} &\mbox{if}& X=C_L, \\
  \mbox{mixed LFSM} &\mbox{if}& X = C_F.
\end{eqnarray*}
\end{definition}

\noindent {\bf Remark.} A given $S \alpha S$ self-similar mixed moving average can be
characterized by different kernels $G$, that is, different integral representations
(\ref{e:mma}). Definition \ref{d:pfsm-cfsm-kernel} and other results of the paper can be
extended to include all $\alpha\in(0,2)$. Their validity depends on the existence of the
so-called ``minimal representations''. If $\alpha\in(1,2)$, there is always at least one
minimal representation (Theorem 4.2 in Pipiras and Taqqu \cite{pipiras:taqqu:2002d}). If
$\alpha\in(0,1]$, a minimal representation exists if some additional conditions are
satisfied (see the Remark on page 436 in Pipiras and Taqqu \cite{pipiras:taqqu:2002d}).
For the sake of simplicity, we chose here not to involve these additional conditions and
hence, we suppose that $\alpha\in(1,2)$, unless otherwise specified. We will not need
here to use explicitly the definition of minimal representations. Minimal representations
are studied, for example, in Rosi{\' n}ski \cite{rosinski:1998} and also used in Pipiras
and Taqqu \cite{pipiras:taqqu:2002d,pipiras:taqqu:2003re}.

\medskip
Mixed LFSMs were studied in Pipiras and Taqqu \cite{pipiras:taqqu:2002s}, Section 7,
where it is shown that they have the following canonical representation.

\begin{proposition}
A $S \alpha S$, $\alpha \in (1,2)$, self-similar mixed moving average $X_\alpha$ is a
mixed LFSM if and only if
\begin{equation}\label{e:mlfsm}
   \left\{
    \begin{array}{ll}
    \int_X\int_\bbR \Big(F_1(x)((t+u)_+^\kappa - u_+^\kappa) + F_2(x)((t+u)_-^\kappa - u_-^\kappa)\Big) M_\alpha(dx,du),& \kappa \ne 0, \\
    \int_X\int_\bbR \Big(F_1(x)\ln\frac{|t+u|}{|u|} + F_2(x)1_{(-t,0)}(u)\Big) M_\alpha(dx,du),& \kappa =
    0,
    \end{array}\right.
\end{equation}
where $F_1,F_2:X\mapsto\bbR$ are some functions and $M_\alpha$ has the control measure
$\nu(dx)du$.
\end{proposition}

If the space $X$ reduces to a single point, then $X_\alpha$ becomes the usual linear
fractional stable motion (LFSM) process if $\kappa\neq 0$, and it becomes a linear
combination of a log-fractional stable motion and a L{\' e}vy stable motion if $\kappa=0$
(see Samorodnitsky and Taqqu \cite{samorodnitsky:taqqu:1994book}, Section 7, for an
introduction to these processes). The process $X_\alpha$ in (\ref{e:mlfsm}) is called a
mixed LFSM because when $\kappa\neq 0$, it differs from a LFSM by the additional variable
$x$.

Our goal is to study integral representations of PFSMs and CFSMs. These processes are
defined in Definition \ref{d:pfsm-cfsm-kernel} in terms of the kernel $G$ in
(\ref{e:mma}) via the sets $C_P$ and $C_L$. An alternative definition of PFSMs and CFSMs
is related to the notion of a flow. A (multiplicative) {\it flow} $\{\psi_c\}_{c>0}$ is a
collection of deterministic maps satisfying
\begin{equation}\label{e:flow0}
\psi_{c_1c_2}(x) = \psi_{c_1}(\psi_{c_2}(x)), \quad \mbox{for all} \ c_1,c_2>0,\ x\in X,
\end{equation}
and $\psi_1(x)=x$, for all $x \in X$. In addition to flows, we shall also use a number of
related real-valued functionals as {\it cocycles}, {\it 1-semi-additive functionals} and
{\it 2-semi-additive functionals}. The definitions of these functionals and related
results are given in Section \ref{s:functionals}. We say henceforth that a $S \alpha S$,
$\alpha \in (0,2)$, self-similar process $X_\alpha$ having a mixed moving average
representation (\ref{e:mma}) is {\it generated by a nonsingular measurable flow}
$\{\psi_c\}_{c>0}$ on $(X,{\cal X},\mu)$ (through the kernel function $G$) if, for all
$c>0$,
\begin{equation}\label{e:generated}
c^{-\kappa} G(x,cu) = b_c(x) \left\{{d(\mu\circ \psi_c) \over d\mu}
(x)\right\}^{1/\alpha} G\Big(\psi_c(x), u + g_c(x)\Big) + j_c(x)\quad \mbox{a.e.}\
\mu(dx)du,
\end{equation}
where $\{b_c\}_{c>0}$ is a cocycle for the flow $\{\psi_c\}_{c>0}$ taking values in
$\{-1,1\}$, $\{g_c\}_{c>0}$ is a {\it 1}-semi-additive functional for the flow
$\{\psi_c\}_{c>0}$ and $\{j_c\}_{c>0}$ is a {\it 2}-semi-additive functional for the flow
$\{\psi_c\}_{c>0}$ and the cocycle $\{b_c\}_{c>0}$, and if
\begin{equation}\label{e:full-support}
\mbox{supp}\left\{ G(x,t+u) - G(x,u),t\in\bbR\right\} = X\times\bbR \quad \mbox{a.e.}\
\mu(dx)du.
\end{equation}
This definition can be found in Pipiras and Taqqu \cite{pipiras:taqqu:2003re}. It differs
from that used in Pipiras and Taqqu
\cite{pipiras:taqqu:2002d,pipiras:taqqu:2002s,pipiras:taqqu:2003cy-fr} in the statement
that $\{ j_c \}_{c>0}$ is a {\it 2}-semi-additive functional, but this, it turns out, is
no restriction.

The following is an alternative definition of PFSMs and CFSMs.

\begin{definition}\label{d:pfsm-cfsm-flow}
A $S \alpha S$, $\alpha \in (1,2)$, self-similar mixed moving average $X_\alpha$ is a
PFSM (a CFSM, resp.) if its minimal representation is generated by a periodic (cyclic,
resp.) flow.
\end{definition}

We now introduce a number of concepts related to Definition \ref{d:pfsm-cfsm-flow}. For
more details, see Pipiras and Taqqu \cite{pipiras:taqqu:2003cy-fr}. A measurable flow
$\{\psi_c\}_{c>0}$ on $(X,{\cal X},\mu)$ is called {\it periodic} if $X=P$ $\mu$-a.e.\
where $P$ is the set of periodic points of the flow defined as
\begin{equation}\label{e:P}
    P =  \{x:\exists\ c = c(x) \neq 1: \psi_c(x) = x\}.
\end{equation}
It is called {\it cyclic} if $X=L$ $\mu$-a.e.\ where $L$ is the set of cyclic points of
the flow defined by
\begin{equation}\label{e:L}
    L = P\setminus F,
\end{equation}
where
\begin{equation}\label{e:F}
    F = \{x: \psi_c(x) = x\ \mbox{for all}\ c>0\}
\end{equation}
is the set of the fixed points.

Note that the sets $C_L$, $C_P$ and $C_F$ in (\ref{e:cfsm-set}), (\ref{e:pfsm-set}) and
(\ref{e:mlfsm-set}) are defined in terms of the kernel $G$, whereas the sets $L$, $P$ and
$F$ are defined in terms of the flow. If the representation is minimal, one has $C_L=L$,
$C_P=P$ and $C_F=F$ $\mu$-a.e.\ (Theorem 3.2, Pipiras and Taqqu
\cite{pipiras:taqqu:2003cy-fr}).

A cyclic flow can also be characterized as a flow which is null-isomorphic (mod 0) to the
flow
\begin{equation}\label{e:cyclic-flow-repres}
\widetilde{\psi}_c(z,v) = (z,\{v + \ln c\}_{q(z)})
\end{equation}
on $(Z\times [0,q(\cdot)),{\cal Z}\times {\cal B}([0,q(\cdot))),\sigma(dz)dv)$, where
$q(z)>0$ a.e.\ is a measurable function. Null-isomorphic (mod $0$) means that there are
two null sets $N\subset X$ and $\widetilde N\subset Z\times [0,q(\cdot))$, and a Borel
measurable, one-to-one, onto and nonsingular map with a measurable nonsingular inverse (a
so-called ``null-isomorphism'') $\Phi:Z\times [0,q(\cdot))\setminus \widetilde N\mapsto X
\setminus N$ such that
\begin{equation}
\psi_c(\Phi(z,v)) = \Phi(\widetilde{\psi}_c(z,v)) \label{e:isomorphic}
\end{equation}
for all $c>0$ and $(z,v)\in Z\times [0,q(\cdot))\setminus \widetilde N$. The null sets
$N$ and $\widetilde N$ are required to be invariant under the flows $\psi_c$ and
$\widetilde{\psi}_c$, respectively. This result, established in Theorem 2.1 of Pipiras
and Taqqu \cite{pipiras:taqqu:2003cy}, will be used in the sequel to establish the
canonical representation (\ref{e:cpfsm-G}) for PFSMs.


\section{Canonical representations for PFSM}
\label{s:canonical-pfsm}

We show in this section that a PFSM can be characterized as a self-similar mixed moving
average represented in one of the explicit ways specified below. We say that these
representations are {\it canonical} for a PFSM. Canonical representations for a PFSM
allow the construction of specific examples of PFSMs and also help to identify them. The
first canonical representation is given in the next result.

\begin{theorem}\label{t:identification-pfsm-canonical}
A $S\alpha S$, $\alpha\in (1,2)$, self-similar mixed moving average $X_\alpha$ is a PFSM
if and only if $X_\alpha$ can be represented by the sum of two independent processes:

\smallskip
$(i)$ The first process has the representation
$$
\int_Z\int_{[0,q(z))}\int_\bbR \Big\{ \Big( b_1(z)^{[v + \ln |t+u|]_{q(z)}}
F_1(z,\{v+\ln|t+u|\}_{q(z)})\, (t+u)_+^\kappa \hspace{2in}
$$
$$
\hspace{2in} -\,  b_1(z)^{[v + \ln |u|]_{q(z)}} F_1(z,\{v+\ln|u|\}_{q(z)})\,
u_+^\kappa\Big)
$$
$$
\hspace{.5in} +\, \Big( b_1(z)^{[v + \ln |t+u|]_{q(z)}} F_2(z,\{v+\ln|t+u|\}_{q(z)})\,
(t+u)_-^\kappa \hspace{2in}
$$
$$
\hspace{2in} -\, b_1(z)^{[v + \ln |u|]_{q(z)}} F_2(z,\{v+\ln|u|\}_{q(z)})\,
u_-^\kappa\Big)
$$
\begin{equation}\label{e:representation}
\hspace{0in} +\, 1_{\{b_1(z) = 1\}} 1_{\{\kappa = 0\}} F_3(z) \ln\frac{|t+u|}{|u|}
\Big\}\, M_\alpha(dz,dv,du),
\end{equation}
where $(Z,{\cal Z},\sigma)$ is a standard Lebesgue space, $b_1(z)\in\{-1,1\}$, $q(z)>0$
a.e.\ $\sigma(dz)$ and $F_1,F_2:Z\times [0,q(\cdot))\mapsto \bbR$, $F_3:Z\mapsto \bbR$
are measurable functions, and $M_\alpha$ has the control measure $\sigma(dz)dvdu$.

\smallskip
$(ii)$ The second process has the representation
\begin{equation}\label{e:mlfsm2}
   \left\{
    \begin{array}{ll}
    \int_Y\int_\bbR \Big(F_1(y)((t+u)_+^\kappa - u_+^\kappa) + F_2(y)((t+u)_-^\kappa - u_-^\kappa)\Big) M_\alpha(dy,du),& \kappa \ne 0, \\
    \int_Y\int_\bbR \Big(F_1(y)\ln\frac{|t+u|}{|u|} + F_2(y)1_{(-t,0)}(u)\Big) M_\alpha(dy,du),& \kappa =
    0,
    \end{array}\right.
\end{equation}
where $(Y,{\cal Y},\nu)$ is a standard Lebesgue space, $F_1,F_2:Y\mapsto\bbR$ are some
functions and $M_\alpha$ has the control measure $\nu(dy)du$.
\end{theorem}

Observe that the processes (\ref{e:mlfsm2}) are the mixed LFSM (\ref{e:mlfsm}) introduced
in Section \ref{s:pfsm-cfsm}. As we will see below (Corollary
\ref{c:canonical-periodic-alt}), Theorem \ref{t:identification-pfsm-canonical} is also
true without part $(ii)$. It is convenient, however, to state it with Part $(ii)$ because
this facilitates the identification of PFSMs and helps in understanding the distinction
between PFSMs and CFSMs. The proof of Theorem \ref{t:identification-pfsm-canonical} can
be found in Section \ref{s:proofs-main}. It is based on results on flow functionals
established in Section \ref{s:represen-cyclic} and on the following proposition.

\begin{proposition}\label{p:representation-cyclic-flow}
If $X_\alpha$ is a $S\alpha S$, $\alpha\in(0,2)$, self-similar mixed moving average
generated by a cyclic flow, then $X_\alpha$ can be represented by
(\ref{e:representation}).
\end{proposition}

The proof of this proposition is given in Section \ref{s:proofs-main}. It is used in the
proof of Theorem \ref{t:identification-pfsm-canonical} in the following way. If
$X_\alpha$ is a PFSM, it has a minimal representation generated by a periodic flow (see
Definition \ref{d:pfsm-cfsm-flow}). Since the periodic points of a flow consist of cyclic
and fixed points, the process $X_\alpha$ can be expressed as the sum of two processes:
one generated by a cyclic flow and the other generated by an identity flow. Proposition
\ref{p:representation-cyclic-flow} is used to show that the process generated by a cyclic
flow has the representation (\ref{e:representation}). The process generated by an
identity flow has the representation (\ref{e:mlfsm2}) according to Theorem 5.1 in Pipiras
and Taqqu \cite{pipiras:taqqu:2002s}.

\smallskip
The representation (\ref{e:representation}) is not specific to processes generated by
cyclic flows. The next result shows that mixed LFSMs (they are generated by identity
flows) can also be represented by (\ref{e:representation}).

\begin{proposition}\label{p:mlfsm-also-representation}
A mixed LFSM having the representation (\ref{e:mlfsm}) can be represented by
(\ref{e:representation}).
\end{proposition}

\begin{proof}
Consider first the case $\kappa \neq 0$. Taking $b_1(z) \equiv 1$, $q(z)\equiv 1$,
$F_1(z,v) \equiv F_1(z)$ and $F_2(z,v) \equiv F_2(z)$ in (\ref{e:representation}), we
obtain the process
$$
\int_Z\int_0^1 \int_\bbR 1_{[0,1)}(v) \Big(F_1(z)((t+u)_+^\kappa - u_+^\kappa) +
F_2(z)((t+u)_-^\kappa - u_-^\kappa) \Big) M_\alpha(dz,dv,du).
$$
Since the kernel above involves the variable $v$ only through the indicator function
$1_{[0,1)}(v)$ and since the control measure of $M_\alpha(dz,dv,du)$ in variable $v$ is
$dv$, the latter process has the same finite-dimensional distributions as
$$
\int_Z \int_\bbR \Big(F_1(z)((t+u)_+^\kappa - u_+^\kappa) + F_2(z)(t+u)_-^\kappa -
u_-^\kappa) \Big) M_\alpha(dz,du),
$$
which is the representation (\ref{e:mlfsm}) of a mixed LFSM when $\kappa \neq 0$. In the
case $\kappa = 0$, one can arrive at the same conclusion by taking $b_1(z) \equiv 1$,
$q(z)\equiv 1$, $F_1(z,v) = F_2(z)$, $F_2(z,v) = 0$ and $F_3(z) = F_1(z)$. \ \ $\Box$
\end{proof}

\medskip
The next corollary which is an immediate consequence of Theorem
\ref{t:identification-pfsm-canonical} and Proposition \ref{p:mlfsm-also-representation}
states, as indicated earlier, that it is not necessary to include Part $(ii)$ in Theorem
\ref{t:identification-pfsm-canonical}.

\begin{corollary}\label{c:canonical-periodic-alt}
A $S\alpha S$, $\alpha\in (1,2)$, self-similar mixed moving average $X_\alpha$ is a PFSM
if and only if it can be represented by (\ref{e:representation}).
\end{corollary}

\medskip
The following result provides another canonical representation of a PFSM which is often
useful in practice. The difference between this result and Theorem
\ref{t:identification-pfsm-canonical} is that the function $s(z)$ appearing in the
expressions $v+s(z) \ln|u|$ below is not necessarily equal to $1$. One can interpret
$|s(z)|$ as the ``speed'' with which the point $(z,v)$ moves under the multiplicative
flow $\psi_c(z,v) = (z,\{v+s(z)\ln c \}_{q(z)})$. The greater $|s(z)|$, the faster does
the fractional part $\{v+s(z)\ln c \}_{q(z)}$ regenerates itself.

\begin{theorem}\label{t:identification-pfsm-canonical2}
A $S\alpha S$, $\alpha\in (1,2)$, self-similar mixed moving average $X_\alpha$ is a PFSM
if and only if $X_\alpha$ can be represented by
$$
\int_Z\int_{[0,q(z))}\int_\bbR \Big\{ \Big( b_1(z)^{[v + s(z) \ln |t+u|]_{q(z)}}
F_1(z,\{v+s(z) \ln|t+u|\}_{q(z)})\, (t+u)_+^\kappa \hspace{1.5in}
$$
$$
\hspace{2in} -\,  b_1(z)^{[v + s(z) \ln |u|]_{q(z)}} F_1(z,\{v+s(z)\ln|u|\}_{q(z)})\,
u_+^\kappa\Big)
$$
$$
\hspace{.5in} +\, \Big( b_1(z)^{[v + s(z) \ln |t+u|]_{q(z)}} F_2(z,\{v+s(z)
\ln|t+u|\}_{q(z)})\, (t+u)_-^\kappa \hspace{2in}
$$
$$
\hspace{2in} -\, b_1(z)^{[v + s(z) \ln |u|]_{q(z)}} F_2(z,\{v+ s(z) \ln|u|\}_{q(z)})\,
u_-^\kappa\Big)
$$
\begin{equation}\label{e:representation2}
\hspace{0in} +\, 1_{\{b_1(z) = 1\}} 1_{\{\kappa = 0\}} F_3(z) \ln\frac{|t+u|}{|u|}\Big)
\Big\}\, M_\alpha(dz,dv,du),
\end{equation}
where $(Z,{\cal Z},\sigma)$ is a standard Lebesgue space, $b_1(z)\in\{-1,1\}$, $s(z)\neq
0$, $q(z)>0$ a.e.\ $\sigma(dz)$ and $F_1,F_2:Z\times [0,q(\cdot))\mapsto \bbR$,
$F_3:Z\mapsto \bbR$ are measurable functions, and $M_\alpha$ has the control measure
$\sigma(dz)dvdu$.
\end{theorem}

\begin{proof}
By Corollary \ref{c:canonical-periodic-alt} above, it is enough to show that the process
(\ref{e:representation2}) can be represented by (\ref{e:representation}). As in Example
2.2 of Pipiras and Taqqu \cite{pipiras:taqqu:2003cy}, the flow
$$
\psi_c(z,v) = (z,\{v+s(z)\ln c\}_{q(z)}),\quad c>0,
$$
is cyclic because each point of the space $Z\times[0,q(\cdot))$ comes back to its initial
position in a finite (nonzero) time. Moreover,
$$
b_c(z,v) = b_1(z)^{[v+s(z)\ln c]},\quad c>0,
$$
is a cocycle for the flow $\{\psi_c\}_{c>0}$. Indeed, by using the second relation in
(\ref{e:int-frac-parts}) and the fact that $\{\psi_c\}_{c>0}$ is a multiplicative flow,
we have
\begin{eqnarray*}
  q[v+s(\ln c_1+\ln c_2)]_q &=& v+s(\ln c_1+\ln c_2) - \{v+s(\ln c_1+\ln c_2)\}_q \\
   &=& v+s\ln c_1 - \{v+s\ln c_1\}_q \\
   & & \hspace{.7in} + \{v+s\ln c_1\}_q + v\ln c_2 - \{\{v+s\ln c_1\}_q + s\ln
c_2 \}_q \\
 &=& q [v+s\ln c_1]_q + q [\{v+s\ln c_1\}_q + s\ln c_2]_q.
\end{eqnarray*}
Hence, $$b_1(z)^{q(z)[v+s(z)(\ln c_1+\ln c_2)]_{q(z)}} = b_1(z)^{[v+s(z)\ln c_1]_{q(z)}}
b_1(z)^{[\{v+s(z)\ln c_1\}_{q(z)} + s(z)\ln c_2]_{q(z)}}$$ which shows that $b_c(z,v)$ is
a cocycle. Since $\{b_c\}_{c>0}$ is a cocycle, relation (\ref{e:generated}) holds and the
process (\ref{e:representation2}) is generated by the flow $\{\psi_c\}_{c>0}$. Since the
flow is cyclic, the process (\ref{e:representation2}) has a representation
(\ref{e:representation}) by Proposition \ref{p:representation-cyclic-flow}. \ \ $\Box$
\end{proof}

\bigskip \noindent {\bf Remark.} Observe that the kernel function $G$ corresponding to the
process (\ref{e:representation2}) is defined as
$$
 G(z,v,u) = b_1(z)^{[v + s(z)\ln |u|]_{q(z)}} \Big( F_1(z,\{v+s(z)\ln|u|\}_{q(z)})\,
u_+^\kappa + F_2(z,\{v+s(z)\ln|u|\}_{q(z)})\, u_-^\kappa\Big)
$$
\begin{equation}\label{e:cpfsm-G-s(z)}
+\ 1_{\{b_1(z) = 1\}} 1_{\{\kappa = 0\}} F_3(z) \ln|u|,
\end{equation}
an expression which will be used a number of times in the sequel.


\section{Representations for CFSM}
\label{s:represen-cyclic}

As shown in Pipiras and Taqqu (2003), CFSMs do not have a (nontrivial) mixed LFSM
component, that is, they cannot be expressed as a sum of two independent processes one of
which is a mixed LFSM. Since a mixed LFSM can be represented by (\ref{e:representation})
(Proposition \ref{p:mlfsm-also-representation} above), a CFSM cannot be characterized as
a process having the representation (\ref{e:representation}). The following result can be
used instead.

\begin{proposition}\label{p:determine-cfsm}
Let $\alpha\in (1,2)$.

$(i)$ If $X_\alpha$ is a CFSM, then it can be represented as (\ref{e:representation}) and
also as (\ref{e:representation2}).

$(ii)$ If the process $X_\alpha$ has the representation (\ref{e:representation}) or
(\ref{e:representation2}) and $C_F=\emptyset$ a.e.\ where the set $C_F$ is defined by
(\ref{e:mlfsm-set}) using the representation (\ref{e:representation}) or
(\ref{e:representation2}), then $X_\alpha$ is a CFSM. Moreover, to show that
$C_F=\emptyset$ a.e., it is enough to prove that $(z,0)\notin C_F$ a.e., that is,
$(z,v)\notin C_F$ a.e.\ when $v=0$.
\end{proposition}

\begin{proof}
$(i)$ If $X_\alpha$ is a CFSM, then $X_\alpha$ given by its minimal representation, is
generated by a cyclic flow by Definition \ref{d:pfsm-cfsm-flow}. Hence, $X_\alpha$ has
the representation (\ref{e:representation}) by Proposition
\ref{p:representation-cyclic-flow}.

$(ii)$ If $X_\alpha$ is represented by (\ref{e:representation}), then by Corollary
\ref{c:canonical-periodic-alt}, $X_\alpha$ is a PFSM. By Definition
\ref{d:pfsm-cfsm-kernel}, the PFSM set $C_P$ associated with the representation
(\ref{e:representation}) is the whole space a.e. Since $C_L = C_P\setminus C_F$, the
assumption $C_F=\emptyset$ a.e.\ implies that $C_L$ is the whole space a.e.\ as well.
Therefore, $X_\alpha$ is a CFSM by Definition \ref{d:pfsm-cfsm-kernel}.

The last statement of Part $(ii)$ follows if we show that $(z,v)\in C_F$ if and only if
$(z,0)\in C_F$. Suppose for simplicity that $X_\alpha$ has the representation
(\ref{e:representation}) defined through the kernel $G$ in (\ref{e:cpfsm-G}). Observe
that
\begin{equation}\label{e:Gzero}
G(z,v,u) = e^{-\kappa v} G(z,0,e^vu) - 1_{\{b_1(z)=1\}} 1_{\{\kappa = 0\}} F_3(z)v
\end{equation}
for all $z,v,u$, and that, by making the change of variables $e^v u =w$, $v=v$ and $z=z$
in (\ref{e:Gzero}),
$$
G(z,0,w) = e^{\kappa v} G(z,v,e^{-v}w) + 1_{\{b_1(z)=1\}} 1_{\{\kappa = 0\}} F_3(z)v
e^{\kappa v}
$$
for all $z,v,w$. By using these relations and the definition (\ref{e:mlfsm-set}) of
$C_F$, we obtain that there is $c_n(z,v)\to 1$ ($c_n(z,v)\neq 1$) such that
$$
G(z,v,c_n(z,v) u) = b_n(z,v) G(z,v,u+a_n(z,v)) + d_n(z,v),\quad \mbox{a.e.}\ du,
$$
for some $a_n(z,v),b_n(z,v)\neq 0,d_n(z,v)$, if and only if there is $\widetilde
c_n(z)\to 1$ ($\widetilde c_n(z)\neq 1$) such that
$$
G(z,0,\widetilde c_n(z) u) = \widetilde b_n(z) G(z,0,u+\widetilde a_n(z)) + \widetilde
d_n(z),\quad \mbox{a.e.}\ du,
$$
for some $\widetilde a_n(z),\widetilde b_n(z)\neq 0,\widetilde d_n(z)$. This shows that
$(z,v)\in C_F$ if and only if $(z,0)\in C_F$. \ \ $\Box$
\end{proof}

\medskip

We do not know of any canonical representation for CFSMs, and feel that such
representations may not exist at all.


\section{Equivalent representations and space of integrands}
\label{s:further-repres-and-space}

The following result shows that there are other representations for PFSMs and CFSMs which
are equivalent to (\ref{e:representation}). As we will see below, this result is useful
for comparing PFSMs and CFSMs to other self-similar stable mixed moving averages, and for
characterizing the space of integrands.

\begin{proposition}\label{p:other-represent}
A process $X_\alpha$ has a representation (\ref{e:representation}) if and only if\begin{eqnarray}
  X_\alpha(t) &\stackrel{d}{=} & \int_Z \int_0^{q(z)} \int_\bbR e^{-\kappa v}
                    (K(z,e^v(t+u)) - K(z,e^vu)) M_\alpha(dz,dv,du) \label{e:equivalent1} \\
   &\stackrel{d}{=}& \int_Z \int_1^{e^{q(z)}} \int_\bbR w^{-H}
                    (K(z,w(t+u)) - K(z,w u)) M_\alpha(dz,dw,du) \label{e:equivalent2} \\
   &\stackrel{d}{=}& \int_Z \int_1^{e^{-q(z)}} \int_\bbR w^{H-\frac{2}{\alpha}}
                    (K(z,w^{-1}(t+u)) - K(z,w^{-1} u)) M_\alpha(dz,dw,du) \label{e:equivalent3} \\
   &\stackrel{d}{=}& \int_Z \int_1^{e^{q(z)}} \int_\bbR w^{-H-\frac{1}{\alpha}}
                    (K(z,w t+u)) - K(z,u)) M_\alpha(dz,dw,du) \label{e:equivalent4} \\
   &\stackrel{d}{=}& \int_Z \int_1^{e^{q(z)}} \int_\bbR w^{H-\frac{1}{\alpha}}
                    (K(z,w^{-1} t+u)) - K(z,u)) M_\alpha(dz,dw,du), \label{e:equivalent5}
\end{eqnarray}
where $M_\alpha(dz,dv,du)$ has the control measure $\sigma(dz)dvdu$, and
$$
K(z,u) = b_1(z)^{[\ln |u|]_{q(z)}} \Big( F_1(z,\{\ln|u|\}_{q(z)})\, u_+^\kappa +
F_2(z,\{\ln|u|\}_{q(z)})\, u_-^\kappa \Big)
$$
\begin{equation}\label{e:cpfsm-G-alt}
+ 1_{\{b_1(z)=1\}} 1_{\{\kappa = 0\}} F_3(z) \ln |u|,
\end{equation}
for some functions $b_1(z)\in\{-1,1\}$, $q(z)>0$ a.e.\ and $F_1,F_2:Z\times [0,q(\cdot))
\mapsto \bbR$, $F_3:Z\mapsto \bbR$.

Moreover, PFSMs and CFSMs can be represented by either of the representation
(\ref{e:equivalent1})--(\ref{e:equivalent5}) with $K$ defined by (\ref{e:cpfsm-G-alt}),
and each of these representations (\ref{e:equivalent1})--(\ref{e:equivalent5}) is
canonical for PFSMs.
\end{proposition}

\begin{proof}
To see that the representations (\ref{e:representation}) and (\ref{e:equivalent1}) are
equivalent, observe that $[v+\ln |u|]_{q(z)} = [\ln |e^vu|]_{q(z)}$, $\{v+\ln
|u|\}_{q(z)} = \{\ln |e^vu|\}_{q(z)}$, and hence
\begin{equation}\label{e:cpfsm-G-G-alt}
G(z,v,u) = e^{-\kappa v} K(z,e^vu),
\end{equation}
where $G$ is the kernel function of (\ref{e:representation}) defined by (\ref{e:cpfsm-G})
and $K$ is defined by (\ref{e:cpfsm-G-alt}). The relations
(\ref{e:equivalent2})--(\ref{e:equivalent5}) follow by making suitable changes of
variables. For example, to obtain (\ref{e:equivalent2}), make the change of variables
$e^v = w$. The last statement of the proposition follows by using Corollary
\ref{c:canonical-periodic-alt} and Proposition \ref{p:determine-cfsm}. \ \ $\Box$
\end{proof}

\medskip
\noindent {\bf Remark.} When the integral over $v$ is $\int_{-\infty}^\infty$ in
(\ref{e:equivalent1}) instead of $\int_0^{q(z)}$, the integrals over $w$ are
$\int_0^\infty$ in (\ref{e:equivalent2})--(\ref{e:equivalent5}) and $K:Z\times\bbR
\mapsto\bbR$ is an arbitrary function (not necessarily of the form
(\ref{e:cpfsm-G-alt})), the resulting processes
(\ref{e:equivalent1})--(\ref{e:equivalent5}) are self-similar stable mixed moving
averages as well. These processes, called {\it dilated fractional stable motions} ({\it
DFSMs}, in short), are generated by the so-called dissipative flows, and are studied in
detail by Pipiras and Taqqu \cite{pipiras:taqqu:2003di}. DFSMs and PFSMs (in particular,
CFSMs) have different finite-dimensional distributions because DFSMs are generated by
dissipative flows and PFSMs are generated by conservative flows (see Theorem 5.3 in
Pipiras and Taqqu \cite{pipiras:taqqu:2002d}). Despite this fact, Proposition
\ref{p:other-represent} shows that the representations of DFSMs and PFSMs have a common
structure.

\medskip
The condition for a PFSM or a CFSM given by (\ref{e:representation}) to be well-defined
is that its kernel function satisfies condition (\ref{e:G_t-in-Lalpha}). By using
Proposition \ref{p:other-represent}, we can replace this condition by one which is often
easier to verify in practice. Let $\alpha\in (0,2)$, $H>0$, $(Z,{\cal Z},\sigma)$ be a
measure space, and $g:Z\mapsto \bbR$ be a function such that $q(z)>0$ a.e.\ $\sigma(dz)$.
Consider the space of functions
\begin{equation}\label{e:space-integrands}
    {\cal C}^q_{\sigma,\alpha,H} = \Big\{ K:Z\times\bbR \mapsto \bbR\ \mbox{such that}\
    \|K\|_{{\cal C}^q_{\sigma,\alpha,H}} <\infty \Big\},
\end{equation}
where
\begin{eqnarray}
  \|K\|^\alpha_{{\cal C}^q_{\sigma,\alpha,H}} &=&
  \int_Z \int_\bbR \int_1^{e^{q(z)}} h^{-\alpha H -1} |K(z,u+h) - K(u)|^\alpha \sigma(dz)dudh
  \label{e:norm-space1}\\
   &=&\int_Z \left( \int_1^{e^{q(z)}} h^{-\alpha H -1} \|\triangle_h K(z,\cdot)\|^\alpha dh \right) \sigma(dz)
  \label{e:norm-space2}
\end{eqnarray}
with $\triangle_h g(\cdot) = g(\cdot + h) - g(\cdot)$.

\begin{proposition}\label{p:well-defined}
A PFSM or CFSM represented by (\ref{e:representation}) is well-defined, that is, the
condition (\ref{e:G_t-in-Lalpha}) holds where $G$ is defined by (\ref{e:cpfsm-G}), if and
only if $K\in {\cal C}^q_{\sigma,\alpha,H}$, where ${\cal C}^q_{\sigma,\alpha,H}$ is
given by (\ref{e:space-integrands}) and $K$ is defined in (\ref{e:cpfsm-G-alt}).
\end{proposition}

\begin{proof}
A process $X_\alpha$ represented by (\ref{e:representation}) is well-defined if and only
if  it is well-defined when represented by (\ref{e:equivalent4}). Since $X_\alpha$ is
self-similar, it is well-defined if and only if the integral
$$
\int_Z \int_1^{e^{q(z)}} \int_\bbR w^{-H-\frac{1}{\alpha}}
                    (K(z,w +u)) - K(z,u)) M_\alpha(dz,dw,du)
$$
corresponding to $X_\alpha(1)$, is well-defined. The latter condition is equivalent to
$K\in {\cal C}^q_{\sigma,\alpha,H}$, where ${\cal C}^q_{\sigma,\alpha,H}$ is given by
(\ref{e:space-integrands}).
 \ \ $\Box$
\end{proof}

\medskip
Though PFSMs and CFSMs were defined for $\alpha\in (1,2)$, Proposition
\ref{p:well-defined} continues to hold for processes represented by
(\ref{e:representation}) when $\alpha \in (0,2)$. For this reason, we defined the space
${\cal C}^q_{\sigma,\alpha,H}$ in (\ref{e:space-integrands}) above for $\alpha \in
(0,2)$. When $Z=\{1\}$ and $\sigma(dz) = \delta_{\{1\}}(dz)$ is the point mass at $z=1$,
we shall use the notation
\begin{equation}\label{e:z=1}
{\cal C}^q_{\alpha,H} := {\cal C}^q_{\sigma,\alpha,H},
\end{equation}
where $q>0$. Thus, $K:\bbR \mapsto \bbR$ is in ${\cal C}^q_{\alpha,H}$ if and only if
\begin{equation}\label{e:space-integrands2}
\int_1^{e^q} h^{-\alpha H -1} \|\triangle_h K(\cdot)\|^\alpha dh < \infty.
\end{equation}
The following result provides sufficient conditions for a function to belong to the space
${\cal C}^q_{\alpha,H}$.

\begin{lemma}\label{l:belong-sufficient}
Let $\alpha\in(0,2)$, $H\in(0,1)$ and $\kappa = H-1/\alpha$, and let $K$ be defined by
(\ref{e:cpfsm-G-alt}) with $Z=\{1\}$ and $\sigma(dz) = \delta_{\{1\}}(dz)$.

\smallskip
$(i)$ Suppose that $\kappa<0$. If $F_1,F_2:[0,q)\mapsto \bbR$ are such that $F_1,F_2$ are
absolutely continuous with derivatives $F_1',F_2'$, and
\begin{equation}\label{e:sufficient1}
    \sup_{u\in[0,q)} |F_i(u)| \leq C, \quad  i=1,2,
\end{equation}
\begin{equation}\label{e:sufficient2}
    \mbox{{\rm ess}} \sup_{u\in[0,q)} |F'_i(u)|\leq C, \quad  i=1,2,
\end{equation}
then $K\in {\cal C}^q_{\alpha,H}$.

\smallskip
$(ii)$ Suppose that $\kappa\geq 0$. If, in addition to $(i)$,
\begin{equation}\label{e:sufficient3}
    F_i(0) = b_1 F_i(q-), \quad i=1,2,
\end{equation}
then $K\in {\cal C}^q_{\alpha,H}$.
\end{lemma}

\begin{proof}
By using (\ref{e:space-integrands2}) and $-\alpha H-1 = -\alpha \kappa -2$, it is enough
to show that
\begin{equation}\label{e:to-show-sufficient}
    \int_1^{e^q} dh\, h^{-\alpha \kappa -2} \int_\bbR du |K(u+h) - K(u)|^\alpha <\infty.
\end{equation}
Since $H\in (0,1)$, we have $\kappa = 0$ only when $\alpha\in (1,2)$. Since, for
$\alpha\in (1,2)$,
$$
\int_\bbR \Big|\ln|u+h| - \ln|u|\Big|^\alpha du = h \int_\bbR \Big|\ln|w+1| -
\ln|w|\Big|^\alpha dw = C h,
$$
where $0<C<\infty$, the function $K(u) = \ln |u|$ satisfies (\ref{e:to-show-sufficient})
when $\kappa =0$. We may therefore suppose that the function $K$ is defined by
(\ref{e:cpfsm-G-alt}) without the last term.

For simplicity, we will prove (\ref{e:to-show-sufficient}) in the case $F_1=F_2=F$, that
is,
\begin{equation}\label{e:to-show-sufficient-K}
    K(u) = b_1^{[\ln |u|]_q} F(\{\ln |u|\}_q) |u|^\kappa,
\end{equation}
where the function $F$ satisfies (\ref{e:sufficient1})--(\ref{e:sufficient2}) in Part
$(i)$ and, in addition, (\ref{e:sufficient3}) in Part $(ii)$. The general case for
arbitrary $F_1$ and $F_2$ can be proved in a similar way. Since $F$ satisfies
(\ref{e:sufficient1}), we have $\int_1^{e^q} dh\, h^{-\alpha \kappa -2} \int_{-N}^N du
|K(u+h) - K(u)|^\alpha <\infty$ for any constant $N$. Since $K(u) = K(-u)$, we have
$\int_1^{e^q} dh\, h^{-\alpha \kappa -2} \int_{N}^\infty du |K(u+h) - K(u)|^\alpha
<\infty$ if and only if $\int_1^{e^q} dh\, h^{-\alpha \kappa -2} \int_{-\infty}^N du
|K(u+h) - K(u)|^\alpha <\infty$. It is therefore enough to show that
\begin{equation}\label{e:to-show-sufficient-alt}
    \int_{N}^\infty |K(u+h) - K(u)|^\alpha du \leq C,\quad h\in (1,e^q),
\end{equation}
where $C$ and $N$ are some constants. Observe that $0<e^{qk}< e^{q(k+1)} - h <
e^{q(k+1)}$ for large enough $k$ and all $h\in [0,q)$. Then, by taking $N=e^{qk_0}$ with
a fixed large $k_0$ in (\ref{e:to-show-sufficient-alt}), we have
$$
\int_{N}^\infty |K(u+h) - K(u)|^\alpha du = \sum_{k=k_0}^\infty
\int_{e^{qk}}^{e^{q(k+1)}-h} |K(u+h) - K(u)|^\alpha du
$$
$$
+ \sum_{k=k_0}^\infty \int_{e^{q(k+1)}-h}^{e^{q(k+1)}} |K(u+h) - K(u)|^\alpha du =: I_1 +
I_2.
$$
We will show that $I_1<\infty$ and $I_2<\infty$.

To show $I_1<\infty$, observe that $qk\leq \ln|u+h| <q(k+1)$ for large $u\in [e^{qk} -
h,e^{q(k+1)} - h)$, and $qk\leq \ln|u| <q(k+1)$ for large $u\in [e^{qk},e^{q(k+1)})$.
Hence, for large $u\in [e^{qk},e^{q(k+1)}-h) \subset [e^{qk} - h,e^{q(k+1)} - h) \cap
[e^{qk},e^{q(k+1)})$, we have
$$
[\ln |u+h| ]_q = [\ln |u| ]_q = k, \quad \{\ln |u+h| \}_q = \ln |u+h| - qk, \quad \{\ln
|u| \}_q = \ln |u| - qk.
$$
By using these relations and since $b_1\in \{-1,1\}$, we obtain that
$$
I_1 =\sum_{k=k_0}^\infty \int_{e^{qk}}^{e^{q(k+1)}-h} \Big| b_1^{[\ln |u+h|]_q}
F(\{\ln|u+h|\}_q) |u+h|^\kappa -  b_1^{[\ln |u|]_q} F(\{\ln |u|\}_q) |u|^\kappa
\Big|^\alpha du
$$
$$
= \sum_{k=k_0}^\infty \int_{e^{qk}}^{e^{q(k+1)}-h} \Big| F(\ln|u+h|-qk) |u+h|^\kappa -
F(\ln |u| - qk) |u|^\kappa \Big|^\alpha du
$$
$$
\leq C \sum_{k=k_0}^\infty \int_{e^{qk}}^{e^{q(k+1)}-h} \Big| F(\ln|u+h|-qk) - F(\ln |u|
- qk) \Big|^\alpha |u+h|^{\kappa\alpha} du
$$
$$
+ C \sum_{k=k_0}^\infty \int_{e^{qk}}^{e^{q(k+1)}-h} |F(\ln |u| - qk)|^\alpha \Big|
|u+h|^\kappa - |u|^\kappa \Big|^\alpha du =: C(I_{1,1} + I_{1,2}).
$$
By using (\ref{e:sufficient2}) and the mean value theorem and by making the change of
variables $u=hw$ below, we obtain that
$$
I_{1,1} \leq C \sum_{k=k_0}^\infty \int_{e^{qk}}^{e^{q(k+1)}-h} \Big|
\ln\frac{|u+h|}{|u|}\Big|^\alpha |u+h|^{\kappa\alpha} du \leq C \int_h^\infty \Big|
\ln\frac{|u+h|}{|u|}\Big|^\alpha |u+h|^{\kappa\alpha} du
$$
$$
= C h^{\kappa \alpha + 1} \int_1^\infty \Big| \ln\frac{|w+1|}{|w|}\Big|^\alpha
|w+1|^{\kappa\alpha} dw \leq C' h^{\kappa \alpha + 1} \int_1^\infty
|w|^{\kappa\alpha-\alpha} dw \leq C,\quad \mbox{for}\ h\in(1,e^q),
$$
since $(\kappa\alpha - \alpha) + 1 = (H-1)\alpha <0$. By using (\ref{e:sufficient1}), we
can similarly show that
$$
I_{1,2} \leq C \int_\bbR \Big||u+h|^\kappa - |u|^\kappa \Big|^\alpha du = C
h^{\kappa\alpha+1} \int_\bbR \Big||w+1|^\kappa - |w|^\kappa \Big|^\alpha dw \leq C',
\quad \mbox{for}\ h\in(1,e^q).
$$
Hence, $I_1\leq C$ for $h\in(1,e^q)$.

We will now show that $I_2<\infty$. Consider first the case $(i)$ where $\kappa<0$. Then,
by using (\ref{e:sufficient1}) and for $h\in (1,e^q)$,
$$
I_2\leq C \sum_{k=k_0}^\infty \int_{e^{q(k+1)}-h}^{e^{q(k+1)}} \Big( |u+h|^{\kappa\alpha}
+ |u|^{\kappa\alpha} \Big) du \leq C' \sum_{k=k_0}^\infty e^{q(k+1)\kappa \alpha}
<\infty,
$$
since $\kappa<0$. Consider now the case $(ii)$ where $\kappa \geq 0$. Observe as above
that, for large $u\in [e^{q(k+1)}-h,e^{q(k+1)})$,
$$
[\ln |u+h| ]_q = k+1,\quad [\ln |u| ]_q = k, \quad \{\ln |u+h| \}_q = \ln |u+h| - q(k+1),
\quad \{\ln |u| \}_q = \ln |u| - qk.
$$
By using these relations and since $b_1\in \{-1,1\}$, we obtain that
$$
I_2 =\sum_{k=k_0}^\infty \int_{e^{q(k+1)}-h}^{e^{q(k+1)}} \Big| b_1^{[\ln |u+h| ]_q}
F(\{\ln|u+h|\}_q) |u+h|^\kappa - b_1^{[\ln |u| ]_q} F(\{\ln |u|\}_q) |u|^\kappa
\Big|^\alpha du
$$
$$
= \sum_{k=k_0}^\infty \int_{e^{q(k+1)}-h}^{e^{q(k+1)}} \Big| b_1 F(\ln|u+h|-q(k+1))
|u+h|^\kappa - F(\ln |u| - qk) |u|^\kappa \Big|^\alpha du
$$
$$
\leq C \sum_{k=k_0}^\infty \int_{e^{q(k+1)}-h}^{e^{q(k+1)}} \Big| b_1 F(\ln|u+h|-q(k+1))
- F(\ln |u| - qk) \Big|^\alpha |u+h|^{\kappa\alpha} du
$$
$$
+ C \sum_{k=k_0}^\infty \int_{e^{q(k+1)}-h}^{e^{q(k+1)}} |F(\ln |u| - qk)|^\alpha \Big|
|u+h|^\kappa - |u|^\kappa \Big|^\alpha du =: C(I_{2,1} + I_{2,2}).
$$
We can show that $I_{2,2}\leq C$ for $h\in (1,e^q)$ as in the case $I_{1,2}$ above.
Consider now the term $I_{2,1}$. By using (\ref{e:sufficient3}), we can extend the
function $F$ to $[q,2q)$ so that $F(u) = b_1 F(u+q)$ for $u\in[0,q)$ or $b_1F(u) =
F(u+q)$ for $u\in[0,q)$ since $b_1 \in \{-1,1\}$, and so that $F$ is absolutely
continuous on $[0,2q)$ with the derivative $F'$ satisfying (\ref{e:sufficient2}). Then,
$$
I_{2,1} = \sum_{k=k_0}^\infty \int_{e^{q(k+1)}-h}^{e^{q(k+1)}} \Big| F(\ln|u+h|-qk) -
F(\ln |u| - qk) \Big|^\alpha |u+h|^{\kappa\alpha} du.
$$
We get $I_{2,1}\leq C$ for $h\in (1,e^q)$ as in the case $I_{1,1}$ above. Hence, $I_2\leq
C$ for $h\in (1,e^q)$ and, since $I_1\leq C$ for $h\in (1,e^q)$ as shown above, we have
$I\leq C$ for $h\in (1,e^q)$ .  \ \ $\Box$
\end{proof}


\section{Examples of PFSMs and CFSMs}
\label{s:examples}

In order for a process given by (\ref{e:representation2}) to be well-defined, its kernel
must belong to the space $L^\alpha(Z\times [0,q(\cdot)) \times\bbR,\sigma(dz)dvdu)$. In
this section, we provide examples of such kernels and hence examples of well-defined
PFSMs. We show that the processes in these examples are also CFSMs.

\begin{example}\label{ex:example-z=1}
Let $\alpha\in(0,2)$ and $H\in(0,1)$ and hence $\kappa = H - 1/\alpha \in
(-1/\alpha,1-1/\alpha)$. The process
\begin{equation}\label{e:example-z=1}
    X_\alpha(t) = \int_0^1 \int_\bbR
    \Big(F(\{v + \ln|t+u| \}_1) (t+u)_+^\kappa - F(\{v + \ln|u| \}_1) u_+^\kappa \Big)
    M_\alpha(dv,du),
\end{equation}
has a representation (\ref{e:representation}) with $Z=\{1\}$,
$\sigma(dz)=\delta_{\{1\}}(dz)$, $b_1(1)=1$, $q(1)=1$, $F_1(1,u)=F(u)$, $F_2(1,u) = 0$
and $F_3(1)=0$. It is well-defined by Lemma \ref{l:belong-sufficient} if the function
$F:[0,1)\mapsto \bbR$ satisfies the conditions (\ref{e:sufficient1}) and
(\ref{e:sufficient2}) when $\kappa<0$ and, in addition, the condition
(\ref{e:sufficient3}) when $\kappa \geq 0$. We can take, for example,
\begin{equation}\label{e:example-kappa<0}
    F(u)= u,\quad u\in[0,1),
\end{equation}
when $\kappa<0$, and
\begin{equation}\label{e:example-any-kappa}
  F(u) = u1_{[0,1/2)}(u) + (1-u) 1_{[1/2,1)}(u),\quad u\in[0,1),
\end{equation}
when no additional conditions on $\kappa$ are imposed. (The function $F$ satisfies
(\ref{e:sufficient3}) because $F(0) = 0 = F(1-)$.) The sufficient conditions of Lemma
\ref{l:belong-sufficient} are not necessary. For example, Lemma 8.1 in Pipiras and Taqqu
\cite{pipiras:taqqu:2002d} shows that the process $X_\alpha$ is also well-defined with
the function
\begin{equation}\label{e:example-kappa<0-2}
    F(u)= 1_{[0,1/2)}(u),\quad u\in[0,1),
\end{equation}
when $\kappa<0$, a function which does not satisfy the (sufficient) conditions of Lemma
\ref{l:belong-sufficient}. When $\alpha \in (1,2)$, the process (\ref{ex:example-z=1})
with the function $F$ in (\ref{e:example-kappa<0}), (\ref{e:example-any-kappa}) or
(\ref{e:example-kappa<0-2}) is therefore a well-defined PFSM.

Note also that, depending on the parameter values $H$ and $\alpha$, the processes
(\ref{e:example-z=1}) may have different sample behavior (even for the same function
$F$). Sample behavior of stable self-similar processes is discussed in Section 12.4 of
Samorodnitsky and Taqqu \cite{samorodnitsky:taqqu:1994book}, and that of general stable
processes is studied in Chapter 10 of that book. When $\kappa>0$, for example, as with
the kernel $F$ in (\ref{e:example-any-kappa}), the process (\ref{e:example-z=1}) is
always sample continuous. When $\kappa<0$, we expect the sample paths to be unbounded
(and hence not continuous) on every interval of positive length for most functions $F$.
For example, when $F$ is c{\` a}dl{\` a}g (that is, right-continuous and with left
limits) on the interval $[0,1]$, then
$$
\sup_{t\in \bbQ} \Big| F(\{v + \ln|t+u| \}_1) (t+u)_+^\kappa - F(\{v + \ln|u| \}_1)
u_+^\kappa\Big| = \infty,
$$
for any $v\in [0,1)$ and $u\in\bbR$ (taking $\bbQ \ni t \to -u$). The unboundedness of
the sample paths follows from Corollary 10.2.4 in Samorodnitsky and Taqqu
\cite{samorodnitsky:taqqu:1994book}. The case $\kappa = 0$ is more difficult to analyze.

\end{example}

\begin{example}\label{ex:example-harmonizable}
Let $\alpha\in(0,2)$ and $H\in (0,1)$. Consider the process
\begin{equation} \label{e:pfsm-connections-harm}
X_\alpha(t) = \int_\bbR \int_0^{2\pi} \int_\bbR \Big( \cos(v+z\ln|t+u|)(t+u)_+^\kappa  -
\cos(v+z\ln|u|)u_+^\kappa  \Big) M_\alpha(dz,dv,du),
\end{equation}
where $M_\alpha(dz,dv,du)$ has the control measure $\lambda(dz)dvdu$. Suppose that the
measure $\lambda(dz)$ is such that
$$
\int_\bbR \lambda(dz)<\infty, \quad \int_\bbR |z|^{\alpha} \lambda(dz)<\infty.
$$
The process (\ref{e:pfsm-connections-harm}) is well-defined since, by the mean value
theorem, one has $|\cos(x) - \cos(z)|\le |x - z|$ and $|\ln|t+u| - \ln|u||\le C|u|^{-1}$
for fixed $t$ and large $|u|$,
$$
\int_\bbR \int_0^{2\pi} \int_\bbR \Big|  \cos(v+z\ln|t+u|)(t+u)_+^\kappa -
\cos(v+z\ln|u|)u_+^\kappa \Big|^\alpha \lambda(dz)dvdu
$$
$$
\le 2^\alpha \int_\bbR \int_0^{2\pi} \int_\bbR \Big| (t+u)_+^\kappa - u_+^\kappa
\Big|^\alpha \lambda(dz)dvdu
$$
$$
+\ 2^\alpha \int_\bbR \int_0^{2\pi} \int_\bbR |u_+|^{\kappa\alpha} \Big|
\cos(v+z\ln|t+u|) - \cos(v+z\ln|u|) \Big|^\alpha \lambda(dz)dvdu
$$
$$
\le 2^\alpha \int_\bbR \lambda(dz) \int_\bbR \Big| (t+u)_+^\kappa - u_+^\kappa
\Big|^\alpha du
$$
$$
+\ 2^\alpha \int_\bbR |z|^\alpha \lambda(dz) \int_\bbR |u_+|^{\kappa\alpha} \Big|
\ln|t+u| - \ln|u| \Big|^\alpha du \ <\infty.
$$
Thus, when $\alpha \in (1,2)$, (\ref{e:pfsm-connections-harm}) is a well-defined PFSM
represented by (\ref{e:representation2}) with $Z=\bbR$, $\sigma(dz) = \lambda(dz)$,
$q(z)=1$, $b_1(z)=1$, $s(z)=z$, $F_1(z,\{u\}_{2\pi}) = \cos (\{u\}_{2\pi}) =\cos (u)$,
$F_2(z,u)=0$ and $F_3(z)=0$.
\end{example}

\medskip
\noindent {\bf Remark.} The PFSM (\ref{e:pfsm-connections-harm}) is related to the
well-known Lamperti transformation and harmonizable processes. Let $\widetilde M(dz,du)$
be a rotationally invariant (isotropic) $S\alpha S$ random measure on $Z\times \bbR$ with
the control measure $\lambda(dz)du$. (Recall that a rotationally invariant $S\alpha S$
random measure $\widetilde M(ds)$ on a space $S$ with the control measure $\widetilde
\nu(ds)$ is a complex-valued random measure satisfying
$$
E\exp\biggl\{ i \Re\Big( \overline{\theta} \int_S f(s) \widetilde M_\alpha(ds) \Big)
\biggr\} = \exp \biggl\{ - |\theta|^\alpha \int_S |f(s)|^\alpha  \widetilde \nu(ds) \Big)
\biggr\}
$$
for all $\theta\in\bbC$ and $f=f_1+if_2\in L^\alpha(S,\widetilde\nu)$. See Samorodnitsky
and Taqqu \cite{samorodnitsky:taqqu:1994book} for more information.) Arguing as in
Example 2.5 of Rosi{\' n}ski (2000), we can show that, up to a multiplicative constant,
the PFSM (\ref{e:pfsm-connections-harm}) has the same finite-dimensional distributions as
the real part of the complex-valued process
\begin{equation}
\int_\bbR \int_\bbR \Big(  e^{iz\ln|t+u|} (t+u)_+^\kappa- e^{iz\ln|u|}u_+^\kappa  \Big)
\widetilde M_\alpha(dz,du). \label{e:pfsm-connections-harm-C}
\end{equation}

The process (\ref{e:pfsm-connections-harm-C}), can also be constructed by the following
procedure. Consider the so-called ``harmonizable'' {\it stationary} process
$$
Y_\alpha^1(t) = \int_\bbR e^{izt} \widetilde M(dz),
$$
where $\widetilde M(dz)$ is rotationally invariant and has the control measure $\mu(dz)$.
By applying the Lamperti transformation to the stationary process $Y_\alpha^1$ (see
Samorodnitsky and Taqqu \cite{samorodnitsky:taqqu:1994book}), one obtains a self-similar
process
$$
Y_\alpha^2(t) = \int_\bbR  e^{iz\ln t} t^H\widetilde M(dz),\ t>0.
$$
The self-similar process $Y_\alpha^2$ can be made also stationary (in the sense of
generalized processes) by introducing an additional variable in its integral
representation, namely,
$$
Y_\alpha^3(t) = \int_\bbR\int_\bbR  e^{iz\ln |t+u|}(t+u)_+^\kappa \widetilde M(dz,du),
$$
where the measure $\widetilde M(dz,du)$ is the same as in
(\ref{e:pfsm-connections-harm-C}). The stationary increments PFSM process
(\ref{e:pfsm-connections-harm-C}) can be obtained from the stationary process
$Y_\alpha^3$ in the usual way by considering $Y_\alpha^3(t) - Y_\alpha^3(0)$.

\medskip

\begin{proposition}\label{p:pfsm-are-cfsp}
The PFSMs of Example \ref{ex:example-z=1} (with $F$ given by either of
(\ref{e:example-kappa<0}), (\ref{e:example-any-kappa}) or (\ref{e:example-kappa<0-2}))
and Example \ref{ex:example-harmonizable} are CFSMs.
\end{proposition}

\begin{proof}
Consider first the process (\ref{e:example-z=1}) defined through the kernel function
$$
G(v,u) = F(\{v+\ln |u|\}_1)u_+^\kappa,
$$
where $F$ is given by either of (\ref{e:example-kappa<0}), (\ref{e:example-any-kappa}) or
(\ref{e:example-kappa<0-2}). By Proposition \ref{p:determine-cfsm}, $(ii)$, it is enough
to show that $0\notin C_F$ where the set $C_F$ is defined by (\ref{e:mlfsm-set}). If
$0\in C_F$, then there is $c_n\to 1$ ($c_n\neq 1$) such that
\begin{equation}\label{e:cannot}
F(\{\ln |c_nu|\}_1) u_+^\kappa = b_n F(\{\ln |u+a_n|\}_1) (u+a_n)_+^\kappa + d_n\quad
\mbox{a.e.}\ du,
\end{equation}
for some $a_n,b_n\neq 0,d_n$. By taking large enough negative $u$ such that $u_+^\kappa =
0$ and $(u+a_n)_+^\kappa=0$, we get that $d_n=0$ and hence (\ref{e:cannot}) becomes
\begin{equation}\label{e:cannot2}
F(\{\ln |c_n u|\}_1) u_+^\kappa = b_n F(\{\ln |u+ a_n|\}_1)(u+a_n)_+^\kappa,\quad
\mbox{a.e.}\ du.
\end{equation}
We shall distinguish between the cases $\kappa<0$ and $\kappa\geq 0$. Observe that we
need to consider the functions (\ref{e:example-kappa<0}), (\ref{e:example-any-kappa}) and
(\ref{e:example-kappa<0-2}) when $\kappa<0$, and only the function
(\ref{e:example-any-kappa}) when $\kappa\geq 0$.

If $\kappa<0$, since $F$ is bounded and not identically zero, by letting $u\to 0$ in
(\ref{e:cannot2}), we obtain that $a_n=0$. Hence, when $\kappa<0$, $F(\{\ln |c_n
u|\}_1)1_{(0,\infty)}(u) = b_n F(\{\ln |u|\}_1)1_{(0,\infty)}(u)$ a.e.\ $du$ or, by
setting $e_n = \ln c_n$ and $v = \ln u$,
\begin{equation}\label{e:cannot3}
    F(\{e_n+v\}_1) = b_nF(\{v\}_1),\quad \mbox{a.e.}\ dv,
\end{equation}
for $e_n\to 0$ ($e_n\neq 0$) and $b_n\neq 0$. Neither of the functions $F$ in
(\ref{e:example-kappa<0}), (\ref{e:example-any-kappa}) or (\ref{e:example-kappa<0-2})
satisfies the relation (\ref{e:cannot3}). For example, if the function
(\ref{e:example-kappa<0}) satisfies (\ref{e:cannot3}), then $e_n+v = b_n v$ for $v\in
[0,1-e_n)$ and some $0<e_n<1$ which is a contradiction (for example, if $v=0$, we get
$e_n=0$). Indeed, if $\kappa\geq 0$, we can also get $a_n=0$ in (\ref{e:cannot2}). If,
for example, $a_n<0$, then since $(u+a_n)_+^\kappa=0$ for $u\in [0,-a_n)$, we have
$F(\{\ln |c_n u|\}_1) = 0$ for $u\in [0,-a_n)$ or that the function $F(v)=0$ on an
interval of $[0,1)$. The function (\ref{e:example-any-kappa}) does not have this
property.

\medskip
Consider now the process (\ref{e:pfsm-connections-harm}) defined through the kernel
function
$$
G(z,v,u) = \cos(v+z\ln|u|) u_+^\kappa.
$$
By Proposition \ref{p:determine-cfsm}, $(ii)$, it is enough to show that $(z,0)\notin
C_F$ a.e.\ $dz$. If $(z,0)\in C_F$, then there is $c_n=c_n(z)\to 1$ ($c_n\neq 1$) such
that
$$
\cos(z\ln|c_nu|) u_+^\kappa = b_n \cos(z\ln|u+a_n|) (u+a_n)_+^\kappa + d_n,\quad
\mbox{a.e.}\ du,
$$
for some $a_n=a_n(z)$, $b_n=b_n(z)\neq 0$, $d_n=d_n(z)$. If $z\neq 0$, we may argue as
above to obtain $d_n=0$ and $a_n=0$. Then,
$$
    \cos(ze_n+v) = b_n \cos(v),\quad \mbox{a.e.}\ dv,
$$
for $e_n\to 0$ ($e_n\neq 0$) and $b_n\neq 0$. This relation cannot hold when $z\neq 0$,
since $e_n\neq 0$, $e_n\to 0$ and because of the shape of the function $\cos(v)$. \ \
$\Box$
\end{proof}

\bigskip \noindent{\bf Remark.} Observe that the $S \alpha S$ $H$--self-similar processes of
Examples \ref{ex:example-z=1} and \ref{ex:example-harmonizable} are well-defined when
$\alpha\in (0,2)$ and $H\in (0,1)$. By Corollary 7.1.1 in Samorodnitsky and Taqqu
\cite{samorodnitsky:taqqu:1994book}, self-similar stable processes can also be defined
when $\alpha\in (0,1)$ and $1\leq H\leq 1/\alpha$, and when $\alpha\in [1,2)$ and $H=1$.
We do not know of examples of processes having a representation (\ref{e:representation})
for these ranges of $\alpha$ and $H$.


\section{Uniqueness results for PFSMs and CFSMs}
\label{s:unique}

We are interested in determining whether two given PFSMs or CFSMs are in fact the same
process. Since we don't want to distinguish between processes which differ by a
multiplicative constant, we will say that two processes $X(t)$ and $Y(t)$ are {\it
essentially identical} if $X(t)$ and $cY(t)$ have the same finite-dimensional
distributions for some constant $c$. If the processes are not essentially identical, we
will say that they are {\it essentially different}.

The next result can often be used to conclude that two PFSMs and CFSMs are essentially
different.

\begin{theorem}\label{t:uniqueness-general}
Suppose that $X_\alpha$ and $\widetilde X_\alpha$ are two PFSMs or CFSMs having
representations (\ref{e:representation2}): the process $X_\alpha$ on the space $Z\times
[0,q(\cdot))\times \bbR$ with the kernel function $G$ defined through the functions
$b_1,s,q,F_1,F_2,F_3$, and the process $\widetilde X_\alpha$ on the space $\widetilde
Z\times [0,\widetilde q(\cdot))\times \bbR$ with the kernel $\widetilde G$ define through
the functions $\widetilde b_1,\widetilde s,\widetilde q,\widetilde F_1,\widetilde
F_2,\widetilde F_3$. If $X_\alpha$ and $\widetilde X_\alpha$ are essentially identical,
then there are maps $h:Z\mapsto \bbR\setminus \{0\}$, $\psi:Z\mapsto \widetilde Z$,
$k:Z\mapsto (0,\infty)$, $g,j:Z\mapsto \bbR$ such that
\begin{equation}\label{e:uniqueness-general}
G(z,0,u) = h(z) \widetilde G(\psi(z),0,k(z)u+g(z)) + j(z),\quad \mbox{a.e.}\
\sigma(dz)du.
\end{equation}
\end{theorem}

\noindent {\bf Remark.} The use of $v=0$ in (\ref{e:uniqueness-general}) should not be
surprising because the function $G(z,v,u)$ in (\ref{e:cpfsm-G-s(z)}) can be expressed
through $G(z,0,w)$. Indeed, as in (\ref{e:Gzero}), it follows from (\ref{e:cpfsm-G-s(z)})
that
\begin{equation}\label{e:G-G-v=0}
    G(z,v,u) = e^{-\kappa v/s(z)} G(z,0,e^{v/s(z)}u) - 1_{\{b_1(z) =1\}} 1_{\{\kappa = 0\}} F_3(z) v/s(z),
    \quad \mbox{for all}\ z,v,u.
\end{equation}

\smallskip
\begin{proof}
Let $G$ and $\widetilde G$ be the kernel functions for the processes $X_\alpha$ and
$\widetilde X_\alpha$, respectively, as defined in the theorem. By Theorem 5.2 in Pipiras
and Taqqu (2002), there are maps $\psi=(\psi_1,\psi_2):Z\times [0,q(\cdot)) \mapsto
\widetilde Z\times [0,\widetilde q(\cdot))$, $h:Z\times [0,q(\cdot)) \mapsto
\bbR\setminus \{0\}$ and $g,j:Z\times [0,q(\cdot)) \mapsto \bbR$ such that
\begin{equation}\label{e:G-Gpsi}
G(z,v,u) = h(z,v) \widetilde G(\psi_1(z,v),\psi_2(z,v),u+g(z,v))+j(z,v)
\end{equation}
a.e.\ $\sigma(dz)dvdu$. By applying (\ref{e:G-G-v=0}) to both sides of (\ref{e:G-Gpsi})
and replacing $u$ by $e^{-v/s(z)}u$, we get
$$
G(z,0,u) = e^{\kappa v/s(z)} e^{-\kappa \psi_2(z,v)/\widetilde s(\psi_1(z,v))} h(z,v)
\widetilde G(\psi_1(z,v),0,e^{\psi_2(z,v)/\widetilde s(\psi_1(z,v))}
(e^{-v/s(z)}u+g(z,v)))
$$
$$
- e^{\kappa v/s(z)} h(z,v) 1_{\{\widetilde b_1(\psi_1(z,v)) = 1\}} 1_{\{\kappa=0\}}
\widetilde F_3(\psi_1(z,v)) \psi_2(z,v)/\widetilde s(\psi_1(z,v))
$$
$$
+ e^{\kappa v/s(z)} 1_{\{b_1(z)=1\}} 1_{\{\kappa =0\}} F_3(z)v/s(z) + e^{\kappa v/s(z)}
j(z,v)
$$
a.e.\ $\sigma(dz)dvdu$. By making the change of variables $e^vu = w$, $z=z$, $v=v$, we
obtain that the last relation holds a.e.\ $\sigma(dz)dvdw$. By fixing $v$, for which this
relation holds a.e.\ $\sigma(dz)dw$, we obtain (\ref{e:uniqueness-general}).
 \ \ $\Box$
\end{proof}

\medskip
The following result shows that we can obtain an ``if and only if'' condition in Theorem
\ref{t:uniqueness-general} in the case when $Z = \{1\}$, $\widetilde Z = \{1\}$, $s =
\widetilde s = 1$ and $q = \widetilde q$. More general cases, for example, even when $Z =
\{1\}$, $\widetilde Z = \{1\}$, $s = \widetilde s = 1$ but $q \neq \widetilde q$, are
much more difficult to analyze. The measure $\delta_{\{1\}}$ below denotes the point mass
at the point $\{1\}$.

\begin{corollary}\label{c:uniqueness-special}
Suppose the conditions of Theorem \ref{t:uniqueness-general} above hold with $Z = \{1\}$,
$\sigma = \delta_{\{1\}}$, $\widetilde Z = \{1\}$, $\widetilde \sigma = \delta_{\{1\}}$,
$s = \widetilde s = 1$ and $q=\widetilde q$. Then, the processes $X_\alpha$ and
$\widetilde X_\alpha$ are essentially identical if and only if there are constants $h\neq
0$, $k>0$, $g,j\in\bbR$ such that
\begin{equation}\label{e:uniqueness-special}
G(0,u) = h \widetilde G(0,k\, u+g) + j,\quad \mbox{a.e.}\ du,
\end{equation}
or equivalently
$$
b_1^{[\ln |u|]_{q}} \Big( F_1(\{\ln|u|\}_{q})\, u_+^\kappa + F_2(\{\ln|u|\}_{q})\,
u_-^\kappa\Big) +\ 1_{\{b_1 = 1\}} 1_{\{\kappa = 0\}} F_3 \ln|u|
$$
$$
= h \widetilde b_1^{[\ln |ku+g|]_{\widetilde q}} \Big( \widetilde
F_1(\{\ln|ku+g|\}_{q})\, (ku+g)_+^\kappa + \widetilde F_2(\{\ln|ku+g|\}_{q})\,
(ku+g)_-^\kappa\Big)
$$
\begin{equation}\label{e:uniqueness-special2}
+\ 1_{\{\widetilde b_1 = 1\}} 1_{\{\kappa = 0\}} h \widetilde F_3 \ln|ku+g| + j,\quad
\mbox{a.e.}\ du.
\end{equation}
Moreover, the processes $X_\alpha$ and $\widetilde X_\alpha$ have identical
finite-dimensional distributions if and only if $|h|=1$.
\end{corollary}

\begin{proof}
The ``only if'' part follows from Theorem \ref{t:uniqueness-general}. We now show the
``if'' part. By using (\ref{e:G-G-v=0}), we have, for $v\in [0,q)$,
$$
G(v,u) = e^{-\kappa v} G(0,e^{v}u) - 1_{\{b_1 =1\}} 1_{\{\kappa = 0\}} F_3 v
$$
$$
= h e^{-\kappa v} \widetilde G(0,ke^vu+g) + e^{-\kappa v} j - 1_{\{b_1 =1\}} 1_{\{\kappa
= 0\}} F_3 v
$$
$$
= h \Big( e^{-\kappa v} \widetilde G(0,ke^vu+g) - 1_{\{\widetilde b_1 =1\}} 1_{\{\kappa =
0\}} \widetilde F_3 v \Big)
$$
$$
+ h 1_{\{\widetilde b_1 =1\}} 1_{\{\kappa = 0\}} \widetilde F_3 v - 1_{\{b_1 =1\}}
1_{\{\kappa = 0\}} F_3 v
$$
$$
= h\widetilde G(v,ku+e^{-v}g) + \widetilde F(v) = h\widetilde G(\{v+\ln
k\}_q,u+e^{-v}k^{-1}g) + \widetilde F(v)
$$
a.e.\ $du$, for some function $\widetilde F(v)$. Hence, by making the change of variables
$u+e^{-v}k^{-1}g \to u$ and $\{v+\ln k\}_q\to v$ below,
$$
X_\alpha(t) \stackrel{d}{=} \int_0^q \int_\bbR \Big( G(v,t+u) - G(v,u) \Big)
M_\alpha(dv,du)
$$
$$
= h\int_0^q \int_\bbR \Big( \widetilde G(\{v+\ln k\}_q,t+u+e^{-v}k^{-1}g) - \widetilde
G(\{v+\ln k\}_q,u+e^{-v}k^{-1}g) \Big) M_\alpha(dv,du)
$$
$$
\stackrel{d}{=} h \int_0^q \int_\bbR \Big( \widetilde G(v,t+u) - \widetilde G(v,u) \Big)
M_\alpha(dv,du) \stackrel{d}{=} h \widetilde X_\alpha(t).
$$
This relation also implies the last statement of the result. \ \ $\Box$
\end{proof}

\bigskip
We now apply Theorem \ref{t:uniqueness-general} and Corollary \ref{c:uniqueness-special}
to examples of PFSMs in Section \ref{s:examples}.

\begin{proposition}\label{p:different}
The four PFSMs considered in Examples \ref{ex:example-z=1} and
\ref{ex:example-harmonizable} are essentially different.
\end{proposition}

\begin{proof}
Let $X_\alpha$ and $Y_\alpha$ be two PFSMs of Example \ref{ex:example-z=1} defined
through two different functions $F_1$ and $F_2$ in (\ref{e:example-kappa<0}),
(\ref{e:example-any-kappa}) or (\ref{e:example-kappa<0-2}). To show that $X_\alpha$ and
$Y_\alpha$ are essentially different, we can suppose that $\kappa<0$ because only the
function (\ref{e:example-any-kappa}) involved $\kappa\geq 0$. If $X_\alpha$ and
$Y_\alpha$ are essentially identical, it follows from Corollary
\ref{c:uniqueness-special} that
$$
F_1(\{\ln|u|\}_{q})\, u_+^\kappa = h  F_2(\{\ln|ku+g|\}_{q})\, (ku+g)_+^\kappa + j
$$
a.e.\ $du$, for $h_1\neq 0$, $k>0$, $g,j\in\bbR$. By arguing as in the proof of
Proposition \ref{p:pfsm-are-cfsp}, we can obtain that $j=0$ and $g=0$. Then, after the
change of variables $\ln u = v$, we have
$$
F_1(\{v\}_{q}) = h_1  F_2(\{k_1 + v\}_{q})
$$
a.e.\ $dv$, for some $h_1\neq 0$, $k_1\neq 0$. This relation does not hold for any two
different functions $F_1$ and $F_2$ in (\ref{e:example-kappa<0}),
(\ref{e:example-any-kappa}) or (\ref{e:example-kappa<0-2}).

Suppose now that $X_\alpha$ and $Y_\alpha$ are the PFSMs of Examples
\ref{ex:example-harmonizable} and \ref{ex:example-z=1}, respectively, defined through the
kernel functions $G_1(z,v,u) = \cos(v+\ln|u|)u_+^\kappa$ and $G_2(z,v) = F(\{v + \ln
|u|\})u_+^\kappa$, where $F$ is defined by (\ref{e:example-kappa<0}),
(\ref{e:example-any-kappa}) or (\ref{e:example-kappa<0-2}). If $X_\alpha$ and $Y_\alpha$
are essentially identical, it follows from Theorem \ref{t:uniqueness-general} that
$$
\cos(z\ln |u|) u_+^\kappa = k(z)F(\{\ln |k(z)u + g(z)|\}_1)(k(z)u+g(z))_+^\kappa + j(z)
$$
a.e.\ $du$, for some $h(z)\neq 0$, $k(z)>0$ and $g(z),j(z)\in\bbR$. When $z\neq 0$, by
arguing as in the proof of Proposition \ref{p:pfsm-are-cfsp}, we get that $j(z)=0$ and
$g(z)=0$. Then, by making the change of variables $\ln u = v$, we have
$$
\cos(zv) = h_1(z) F(\{k_1(z) + v\}_1)
$$
a.e.\ $dv$, for some $h_1(z)\neq 0$, $k_1(z)\neq 0$. The function $F$ in
(\ref{e:example-kappa<0}), (\ref{e:example-any-kappa}) or (\ref{e:example-kappa<0-2})
does not satisfy this relation. Hence, $X_\alpha$ and $Y_\alpha$ are essentially
different.
 \ \ $\Box$
\end{proof}


\section{Functionals of cyclic flows} \label{s:functionals}

We shall characterize here the functionals appearing in (\ref{e:generated}) which are
related to cyclic flows. These results are used in Section \ref{s:proofs-main} below to
establish a canonical representation for PFSMs in Theorem
\ref{t:identification-pfsm-canonical}. We start by providing a precise definition of
flows and related functionals. See Pipiras and Taqqu \cite{pipiras:taqqu:2003re} for
motivation. A flow $\{\psi_c\}_{c>0}$ on a standard Lebesque space $(X,{\cal X},\mu)$ is
a collection of measurable maps $\psi_c:X \to X$ satisfying (\ref{e:flow0}). The flow is
{\it nonsingular} if each map $\psi_c$, $c>0$, is nonsingular, that is, $\mu(A)=0$
implies $\mu(\psi_c^{-1} (A))=0$. It is {\it measurable} if a map $\psi_c(x):(0, \infty)
\times X \to X$ is measurable. A {\it cocycle} $\{b_c\}_{c>0}$ for the flow
$\{\psi_c\}_{c>0}$ taking values in $\{-1,1\}$ is a measurable map $b_c(x):(0, \infty)
\times X \to \{-1,1\}$ such that
\begin{equation}\label{e:cocycle0}
b_{c_1c_2}(x) = b_{c_1}(x) b_{c_2}(\psi_{c_1}(x)), \quad \mbox{for all} \ c_1,c_2>0,\
x\in X.
\end{equation}
A {\it 1-semi-additive functional} $\{g_c\}_{c>0}$ for the flow $\{\psi_c\}_{c>0}$ is a
measurable map $g_c(x):(0, \infty) \times X \to \bbR$ such that
\begin{equation}\label{e:1-semi-add0}
g_{c_1c_2}(x) = c_2^{-1} g_{c_1}(x) + g_{c_2}(\psi_{c_1}(x)), \quad \mbox{for all} \
c_1,c_2>0,\ x\in X.
\end{equation}
A {\it 2-semi-additive functional} $\{j_c\}_{c>0}$ for the flow $\{\psi_c\}_{c>0}$ and a
related cocycle $\{b_c\}_{c>0}$ is a measurable map $j_c(x):(0, \infty) \times X \to
\bbR$ such that
\begin{equation}\label{e:2-semi-add0}
j_{c_1c_2}(x) = c_2^{-\kappa} j_{c_1}(x) + b_{c_1}(x)\left\{\frac{d(\mu\circ
\psi_{c_1})}{d\mu}(x)\right\}^{1/\alpha} j_{c_2}(\psi_{c_1}(x)),\quad \mbox{for all} \
c_1,c_2>0,\ x\in X.
\end{equation}
The Radon-Nikodym derivatives $\widetilde b_c(x)=(d(\mu\circ \psi_c)/d\mu)(x)$ in
(\ref{e:2-semi-add0}) can and will be viewed as a cocycle taking values in $\bbR
\setminus \{0\}$, that is, a measurable map $\widetilde b_c(x):(0, \infty) \times X \to
\bbR \setminus \{0\}$ satisfying $\widetilde b_{c_1c_2}(x) = \widetilde b_{c_1}(x)
\widetilde b_{c_2}(\psi_{c_1}(x))$, for all $c_1,c_2>0$ and $x \in X$ (see Pipiras and
Taqqu \cite{pipiras:taqqu:2003re}).

In the following three lemmas, we characterize cocycles and {\it 1}- and {\it
2}-semi-additive functionals associated with cyclic flows.

\begin{lemma}\label{l:cocycle-cyclic}
Let $\{b_c\}_{c>0}$ be a cocycle taking values in $\{-1,1\}$ for a cyclic flow
$\{\psi_c\}_{c>0}$. Set $\widetilde{b}_c(z,v) = b_c(\Phi(z,v))$ if $(z,v)\in Z\times
[0,q(\cdot))\setminus \widetilde N$, and $\widetilde{b}_c(z,v) =1$ if $(z,v)\in
\widetilde N$, where $\Phi:Z\times [0,q(\cdot))\setminus \widetilde N \mapsto X \setminus
N$ is the null-isomorphism appearing in (\ref{e:isomorphic}). Then, $\{
\widetilde{b}_c\}_{c>0}$ is a cocycle for the cyclic flow $\{\widetilde \psi_c\}_{c>0}$
defined by (\ref{e:cyclic-flow-repres}), and it can be expressed as
\begin{equation} \label{e:cocycle-cyclic-represent}
  \widetilde{b}_c(z,v) = \frac{\widetilde{b}(z,\{v + \ln c\}_{q(z)})}{\widetilde{b}(z,v)}\
  \widetilde b_1(z)^{[v + \ln c]_{q(z)}},
  \end{equation}
for some functions $\widetilde b: Z \times [0,q(\cdot)) \mapsto \{-1,1\}$ and
$\widetilde{b}_1:Z\mapsto \{-1,1\}$.
\end{lemma}

\begin{proof} The result can be deduced from Lemma 2.2 in
Pipiras and Taqqu \cite{pipiras:taqqu:2003cy} by using the following relation between
multiplicative and additive flows: $\{\psi_c\}_{c>0}$ is a multiplicative flow and
$\{b_c\}_{c>0}$ is a related cocycle if and only if $\phi_t:=\psi_{e^t}, t \in \bbR$, is
an additive flow and $a_t:=b_{e^t}, t \in \bbR$, is a related cocycle. (Additive flows
$\{\phi_t\}_{t \in \bbR}$ are such that $\phi_{t_1+t_2}= \phi_{t_1} \circ \phi_{t_2},
t_1, t_2 \in \bbR.$)
 \ \ $\Box$
\end{proof}

\begin{lemma}\label{l:semi-additive-cyclic}
Let $\{g_c\}_{c>0}$ be a {\it 1}-semi-additive functional for a cyclic flow
$\{\psi_c\}_{c>0}$. Set $\widetilde{g}_c(z,v) = g_c(\Phi(z,v))$ if $(z,v)\in Z\times
[0,q(\cdot))\setminus \widetilde N$, and $\widetilde{g}_c(z,v)=0$ if $(z,v)\in \widetilde
N$, where $\Phi:Z\times [0,q(\cdot))\setminus \widetilde N \mapsto X \setminus N$ is the
null-isomorphism appearing in (\ref{e:isomorphic}). Then, $\{\widetilde{g}_c\}_{c>0}$ is
a {\it 1}-semi-additive functional for the cyclic flow $\{\widetilde \psi_c\}_{c>0}$
defined by (\ref{e:cyclic-flow-repres}), and it can be expressed as
\begin{equation} \label{e:1-semi-additive-cyclic-represent}
  \widetilde{g}_c(z,v) = \widetilde{g}(z,\{v + \ln c\}_{q(z)})- c^{-1}  \widetilde{g}(z,v),
  \end{equation}
for some function $\widetilde g:Z\times [0,q(\cdot)) \mapsto \bbR$.
\end{lemma}

\begin{proof} We may suppose without loss of generality $N= \widetilde N = \emptyset$
because $\widetilde{g}_c(z,v)=0$ obviously satisfies the {\it 1}-semi-additive functional
relation (\ref{e:1-semi-add0}) on the set $\widetilde N$ (which is invariant for the flow
$\widetilde \psi_c$). By substituting $x=\Phi(z,v)$ in the equation
(\ref{e:1-semi-add0}), we obtain that
$$
g_{c_1c_2}(\Phi(z,v))=c_2^{-1} g_{c_1} (\Phi(z,v))+ g_{c_2} (\psi_{c_1}(\Phi(z,v)))
$$
and, since $\psi_c \circ \Phi = \Phi \circ \widetilde \psi_c$ by (\ref{e:isomorphic}), we
get
\begin{equation} \label{e:1-semi-additive-tildeg}
\widetilde g_{c_1c_2}(z,v)=c_2^{-1} \widetilde g_{c_1} (z,v)+ \widetilde g_{c_2}
(\widetilde \psi_{c_1}(z,v)).
\end{equation}
Relation (\ref{e:1-semi-additive-tildeg}) shows that $\{ \widetilde g_c\}_{c>0}$ is a
{\it 1}-semi-additive functional for the flow $\{ \widetilde \psi_c\}_{c>0}$. The
expression (\ref{e:1-semi-additive-cyclic-represent}) follows from Proposition 5.1 in
Pipiras and Taqqu \cite{pipiras:taqqu:2003re}. \ \ $\Box$
\end{proof}

\bigskip Let $T$ be an arbitrary index set, e.g.\ $T=(0,\infty)$, and $\{f_t^1\}_{t\in
T}$, $\{f_t^1\}_{t\in T}$ be two collections of measurable functions on $(X,{\cal
X},\mu)$. We will say that $\{f_t^1\}_{t\in T}$ is a version of $\{f_t^2\}_{t\in T}$ if
$\mu \{x:f_t^1(x) \neq f_t^2(x)\} = 0$ for all $t\in T$. We now characterize {\it
2}-semi-additive functionals related to cyclic flows.

\begin{lemma}\label{l:2-semi-additive-cyclic}
Let $\{\widetilde{j}_c\}_{c>0}$ be a {\it 2}-semi-additive functional for a cyclic flow
$\{\psi_c\}_{c>0}$. Set
\begin{equation}\label{e:widetilde-j_c-j_c}
\widetilde{j}_c(z,v) =  \left\{\frac{d(\mu\circ \Phi)}{d(\sigma \otimes \bbL)}(z,v)
\right\}^{1/\alpha}j_c(\Phi(z,v))
\end{equation}
if $(z,v)\in Z\times [0,q(\cdot))\setminus \widetilde N$, and let $\widetilde{j}_c(z,v)$
be defined arbitrarily for $(z,v)\in\widetilde N$, where $\Phi:Z \times [0,q(\cdot))
\setminus \widetilde N \mapsto X \setminus N$ is the isomorphism appearing in
(\ref{e:isomorphic}). Then, $\{\widetilde{j}_c\}_{c>0}$ has a version which is a {\it
2}-semi-additive functional for the flow $\{\widetilde{\psi}_c\}_{c>0}$ defined by
(\ref{e:cyclic-flow-repres}) and the cocycle $\{\widetilde{b}_c\}_{c>0}$ defined in Lemma
\ref{l:cocycle-cyclic}. Moreover, for any $c>0$,
$$
\widetilde{j}_c(z,v) = \widetilde b_c(z,v) \widetilde j(z,\{v + \ln c\}_{q(z)}) -
c^{-\kappa} \widetilde j(z,v)
$$
\begin{equation}\label{e:2-semi-additive-cyclic-repres}
+\, \widetilde j_1(z)(\widetilde{b}(z,v))^{-1} [v + \ln c]_{q(z)} 1_{\{\widetilde
b_1(z)=1\}} 1_{\{\kappa = 0\}}, \quad \mbox{a.e.}\ \sigma(dz)dv,
\end{equation}
where $\widetilde j:Z\times [0,q(\cdot))\mapsto \bbR$, $\widetilde j_1:Z\mapsto \bbR$ are
some measurable functions, and the functions $\widetilde b_1(z)$ and $\widetilde b(z,v)$
appear in (\ref{e:cocycle-cyclic-represent}) of Lemma \ref{l:cocycle-cyclic}.
\end{lemma}

\begin{proof}
We suppose without loss of generality that $N = \widetilde N = \emptyset$. Substituting
$x= \Phi(z,v)$ in the {\it 2}-semi-additive functional equation (\ref{e:2-semi-add0}), we
obtain that
\begin{equation} \label{e:2-semi-additive-j}
j_{c_1c_2}(\Phi(z,v)) = c_2^{-\kappa} j_{c_1}(\Phi(z,v)) + b_{c_1}(\Phi(z,v))
\left\{\frac{d(\mu\circ \psi_{c_1})}{d\mu}(\Phi(z,v)) \right\}^{1/\alpha}
j_{c_2}(\psi_{c_1}(\Phi(z,v))).
\end{equation}
We first show that $\{\widetilde j_c\}_{c>0}$ is an almost {\it 2}-semi-additive
functional, that is, it satisfies relation (\ref{e:2-semi-add0}) for all $c_1,c_2>0$
a.e.\ $\sigma(dz)dv$. Observe that, for any $c>0$ and $\bbL$ denoting the Lebesgue
measure,
\begin{eqnarray}
  \frac{d(\mu\circ \psi_c)}{d\mu}\circ \Phi  &=&  \frac{d((\mu\circ \Phi) \circ \widetilde \psi_c)}
  {d(\mu\circ \Phi)}\nonumber \\
   &=& \frac{d((\sigma \otimes \bbL)\circ \widetilde \psi_c)}{d(\sigma\otimes \bbL)}\,
   \frac{d(\sigma\otimes \bbL)}{d(\mu\circ \Phi)}\,
   \left( \frac{d(\mu\circ \Phi)}{d(\sigma\otimes \bbL)}\circ \widetilde \psi_c\right)  \nonumber\\
   &=& \left( \frac{d(\mu\circ \Phi)}{d(\sigma\otimes \bbL)}\right)^{-1}\left(
   \frac{d(\mu\circ \Phi)}{d(\sigma\otimes \bbL)}\circ \widetilde \psi_c\right), \quad (\sigma
   \otimes \bbL)\mbox{-a.e.}, \label{e:jacobian}
\end{eqnarray}
where in the last equality above we used the identity $d((\sigma \otimes \bbL)\circ
\widetilde \psi_c) / d(\sigma\otimes \bbL) = 1$ $(\sigma \otimes \bbL)$-a.e., which
follows from (\ref{e:cyclic-flow-repres}) because $d\{v+\ln c\}/dv=1$ a.e. By using the
relation (\ref{e:jacobian}), we can write (\ref{e:2-semi-additive-j}) as
$$
\left\{ \frac{d(\mu\circ \Phi)}{d(\sigma \otimes \bbL)}(z,v)\right\}^{1/\alpha}
j_{c_1c_2}(\Phi(z,v)) = c_2^{-\kappa} \left\{ \frac{d(\mu\circ \Phi)}{d(\sigma \otimes
\bbL)}(z,v)\right\}^{1/\alpha} j_{c_1}(\Phi(z,v))
$$
\begin{equation} \label{e:2-semi-additive-j2}
+\ b_{c_1}(\Phi(z,v)) \left\{\frac{d(\mu\circ \Phi)}{d(\sigma \otimes \bbL)}(\widetilde
\psi_c(z,v)) \right\}^{1/\alpha} j_{c_2}(\psi_{c_1}(\Phi(z,v))), \quad \mbox{a.e.} \
\sigma(dz)dv.
\end{equation}
Since $\psi_c\circ \Phi = \Phi \circ \widetilde \psi_c$ by (\ref{e:isomorphic}), $b_c
\circ \Phi = \widetilde b_c$ by using the notation of Lemma \ref{l:cocycle-cyclic} and
$\widetilde j_c = \{d(\mu\circ \Phi)/d(\sigma \otimes \bbL)\}^{1/\alpha} j_c$ by
(\ref{e:widetilde-j_c-j_c}), we deduce from (\ref{e:2-semi-additive-j2}) that, for any
$c_1,c_2>0$,
\begin{equation}\label{e:2-semi-additive-tildej}
\widetilde j_{c_1c_2}(z,v) = c_2^{-\kappa} \widetilde j_{c_1}(z,v) + \widetilde
b_{c_1}(z,v) \widetilde j_{c_2}(\widetilde \psi_{c_1}(z,v)), \quad \mbox{a.e.} \
\sigma(dz)dv,
\end{equation}
that is, $\{\widetilde j_c\}_{c>0}$ is an almost {\it 2}-semi-additive functional. By
Theorem 2.1 in Pipiras and Taqqu \cite{pipiras:taqqu:2003re}, $\{\widetilde j_c\}_{c>0}$
has a version which is a {\it 2}-semi-additive functional.

Since $\{\widetilde j_c\}_{c>0}$ has a version which is a {\it 2}-semi-additive
functional, we may suppose without loss of generality that $\{\widetilde j_c\}_{c>0}$ is
a {\it 2}-semi-additive functional. The expression
(\ref{e:2-semi-additive-cyclic-repres}) for $\{\widetilde j_c\}_{c>0}$ then follows from
Proposition 5.2 in Pipiras and Taqqu \cite{pipiras:taqqu:2003re}.  \ \ $\Box$
\end{proof}


\section{The proofs of Proposition \ref{p:representation-cyclic-flow} and Theorem \ref{t:identification-pfsm-canonical}}
\label{s:proofs-main}

{\sc Proof of Proposition \ref{p:representation-cyclic-flow}:} Suppose that the process
$X_\alpha$ is generated by a cyclic flow $\{\psi_c\}_{c>0}$ on $(X,{\cal X},\mu)$. Then,
by a discussion following Definition \ref{d:pfsm-cfsm-flow}, there are a standard
Lebesgue space $(Z,{\cal Z},\sigma)$, function $q(z)>0$ and a null-isomorphism $\Phi:
Z\times [0,q(\cdot))\mapsto X$ such that
\begin{equation}\label{e:isomorphic-again}
\psi_c(\Phi(z,v)) = \Phi(z,\{v + \ln c\}_{q(z)})
\end{equation}
for all $c>0$ and $(z,v)\in Z\times [0,q(\cdot))$. In other words, the flow
$\{\psi_c\}_{c>0}$ on $(X,\mu)$ is null-isomorphic to the flow $\{\widetilde
\psi_c\}_{c>0}$ on $(Z\times[0,q(\cdot)),\sigma(dz)dv)$ defined by $\widetilde
\psi_c(z,v) = (z,\{v + \ln c\}_{q(z)})$. (We may suppose that the null sets in
(\ref{e:isomorphic}) are empty because, otherwise, we can replace $X$ by $X\setminus N$
in the definition of $X_\alpha$ without changing its distribution.) By replacing $x$ by
$\Phi(z,v)$ in (\ref{e:generated}) and using (\ref{e:isomorphic-again}), we get that, for
all $c>0$,
$$
c^{-\kappa} G(\Phi(z,v), cu) = b_c(\Phi(z,v)) \left\{ {d(\mu \circ \psi_c)\over
d\mu}(\Phi(z,v)) \right\}^{1/\alpha}\times
$$
\begin{equation}
\times\, G \Big(\Phi(\widetilde \psi_c(z,v) ), u + g_c(\Phi(z,v)) \Big) + j_c(\Phi(z,v))
\label{e:for-char-cycle-1}
\end{equation}
a.e.\ $\sigma(dz)dvdu$. By using the relation
\begin{eqnarray}
{d(\mu \circ \psi_c)\over d\mu}\circ \Phi & = &  {d(\mu \circ \Phi  \circ \widetilde
\psi_c)
\over d(\mu\circ \Phi) }             \nonumber \\
& = & \left({d\mu \over d((\sigma\otimes \bbL)\circ \Phi^{-1})} \circ \Phi\circ
\widetilde \psi_c\right)\, {d((\sigma\otimes \bbL)\circ \widetilde \psi_c)
 \over d(\sigma\otimes \bbL)}\,
{d(\sigma\otimes \bbL)
 \over d (\mu \circ \Phi)} \nonumber \\
& = & \left({d\mu \over d((\sigma\otimes \bbL)\circ \Phi^{-1})} \circ \Phi\circ
\widetilde \psi_c \right)\, {d((\sigma\otimes \bbL)\circ \Phi^{-1})
 \over d \mu }\circ \Phi \nonumber  \\
& = & \left({d\mu \over d((\sigma\otimes \bbL)\circ \Phi^{-1})} \circ \Phi\circ
\widetilde \psi_c\right)\, \left( {d \mu \over d((\sigma\otimes \bbL)\circ
\Phi^{-1})}\circ \Phi \right)^{-1},\nonumber
\end{eqnarray}
where $\bbL$ is the Lebesgue measure, setting
\begin{equation}\label{e:widetildeG-G}
\widetilde G(z,v,u) = \left\{\frac{d(\mu\circ \Phi)}{d(\sigma \otimes \bbL)}(z,v)
\right\}^{1/\alpha} G(\Phi(z,v), u)
\end{equation}
and using the notation of Lemmas \ref{l:cocycle-cyclic}, \ref{l:semi-additive-cyclic} and
\ref{l:2-semi-additive-cyclic}, we obtain that, for all $c>0$,
\begin{equation}\label{e:to-solve-for-tilde-G}
c^{-\kappa}\widetilde{G}(z,v,cu) = \widetilde b_c(z,v) \widetilde{G}(\widetilde
\psi_c(z,v), u + \widetilde g_c(z,v)) + \widetilde j_c(z,v)
\end{equation}
a.e.\ $\sigma(dz)dvdu$. We next consider the cases $\kappa\neq 0$ and $\kappa = 0$
separately.

\medskip
{\it The case $\kappa\neq 0$:} By using Lemmas \ref{l:cocycle-cyclic},
\ref{l:semi-additive-cyclic} and \ref{l:2-semi-additive-cyclic}, and setting
\begin{equation}\label{e:hatG}
    \widehat G(z,v,u) = \widetilde b(z,v) \Big( \widetilde G(z,v,u+\widetilde g(z,v)) -
    \widetilde j(z,v)\Big),
\end{equation}
we obtain from (\ref{e:to-solve-for-tilde-G}) that, for all $c>0$,
$$
\widehat G(z,v,cu) = c^\kappa \widetilde b_1(z)^{[v+\ln c]_{q(z)}} \widehat G(z,\{v+\ln
c\}_{q(z)},u)
$$
a.e.\ $\sigma(dz)dvdu$. By making a change of variables $cu = w$, we get that, for all
$c>0$,
$$
\widehat G(z,v,w) = c^\kappa \widetilde b_1(z)^{[v+\ln c]_{q(z)}} \widehat G(z,\{v+\ln
c\}_{q(z)},c^{-1} w)
$$
$$
= |w^{-1}c|^\kappa |w|^\kappa \widetilde b_1(z)^{[v+\ln |w| |w^{-1} c|]_{q(z)}} \widehat
G(z,\{v+\ln |w||w^{-1} c|\}_{q(z)},c^{-1}w)
$$
a.e.\ $\sigma(dz)dvdw$. By the Fubini's theorem, this relation holds a.e.\
$\sigma(dz)dvdwdc$ as well. Then, for $w>0$, by making the change of variables $c = yw$
and then fixing $y=y_0$, we get
$$
\widehat G(z,v,w) = w^\kappa \widetilde b_1(z)^{[v+\ln a_1|w|]_{q(z)}} \widehat
F_1(z,\{v+\ln a_1|w|\}_{q(z)})
$$
a.e.\ $\sigma(dz)dvdw$, for some $a_1>0$ and function $\widehat F_1$. By using identities
$\{v+\ln a_1|w|\}_{q(z)} = \{ \{v+\ln |w|\}_{q(z)} + \ln a_1\}_{q(z)}$ and $[v+\ln
a_1|w|]_{q(z)} = [v+\ln |w|]_{q(z)}+[\{v+\ln |w|\}_{q(z)} + \ln a_1]_{q(z)}$, we can
simplify the last relation as
\begin{equation}\label{e:tilde-F1-not0}
\widehat G(z,v,w) = w^\kappa \widetilde b_1(z)^{[v+\ln |w|]_{q(z)}} F_1(z,\{v+\ln
|w|\}_{q(z)})
\end{equation}
a.e.\ $\sigma(dz)dvdw$, for some function $F_1$. Similarly, for $w<0$, we may get that
\begin{equation} \label{e:tilde-F2-not0}
\widehat G(z,v,w) = w_-^\kappa \widetilde b_1(z)^{[v+\ln |w|]_{q(z)}} F_2(z,\{v+\ln
|w|\}_{q(z)})
\end{equation} a.e.\
$\sigma(dz)dvdw$, for some function $F_2$. Observe now that, by writing characteristic
functions and using (\ref{e:widetildeG-G}) and (\ref{e:hatG}),
$$
\{X_\alpha(t)\}_{t\in\bbR} \stackrel{d}{=} \left\{\int_X\int_\bbR (G(x,t+u) - G(x,u))
M_\alpha(dx,du) \right\}_{t\in\bbR}
$$
$$
\stackrel{d}{=} \left\{\int_{Z}\int_{[0,q(\cdot))} \int_\bbR (\widetilde{G}(z,v,t+u) -
\widetilde{G}(z,v,u)) \widetilde{M}_\alpha(dz,dv,du) \right\}_{t\in\bbR}
$$
\begin{equation}
\stackrel{d}{=} \left\{\int_{Z}\int_{[0,q(\cdot))} \int_\bbR (\widehat{G}(z,v,t+u) -
\widehat{G}(z,v,u)) \widetilde{M}_\alpha(dz,dv,du) \right\}_{t\in\bbR},
\label{e:back-to-G-not0}
\end{equation}
where $\widetilde{M}_\alpha(dz,dv,du)$ is a $S\alpha S$ random measure on
$(Z\times[0,q(\cdot)))\times \bbR$ with the control measure $\sigma(dz)dvdu$. The result
of the theorem when $\kappa\ne 0$ then follows by using (\ref{e:tilde-F1-not0}) and
(\ref{e:tilde-F2-not0}).

\medskip
{\it The case $\kappa=0$:} In this case, by using Lemmas \ref{l:cocycle-cyclic},
\ref{l:semi-additive-cyclic} and \ref{l:2-semi-additive-cyclic}, and the notation
(\ref{e:hatG}), we get that, for all $c>0$,
$$
\widehat G(z,v,cu) = \widetilde b_1(z)^{[v+\ln c]_{q(z)}} \widehat G(z,\{v+\ln
c\}_{q(z)},u) + \widetilde j_1(z) [v+\ln c]_{q(z)} 1_{\{\widetilde b_1(z) = 1\}}
$$
a.e.\ $\sigma(dz)dvdu$. Arguing as in the case $\kappa\neq 0$, we may show that, for
$w>0$,
$$
\widehat G(z,v,w) = \widetilde b_1(z)^{[v+\ln a_1|w|]_{q(z)}} \widehat F_1(z,\{v+\ln
a_1|w|\}_{q(z)}) 1_{(0,\infty)}(w)
$$
$$
+ \widetilde j_1(z) [v+\ln a_1|w|]_{q(z)} 1_{\{\widetilde b_1(z) = 1\}}1_{(0,\infty)}(w),
$$
a.e.\ $\sigma(dz)dvdu$. By using the identities preceding (\ref{e:tilde-F1-not0}), we
conclude that, for $w>0$,
$$
\widehat G(z,v,w) = \widetilde b_1(z)^{[v+\ln |w|]_{q(z)}} F_1(z,\{v+\ln |w|\}_{q(z)})
1_{(0,\infty)}(w)
$$
$$
+ F_3(z) [v+\ln |w|]_{q(z)} 1_{\{\widetilde b_1(z) = 1\}}1_{(0,\infty)}(w),
$$
a.e.\ $\sigma(dz)dvdu$, for some functions $F_1$ and $F_3$. Similarly, for $w<0$,
$$
\widehat G(z,v,w) = \widetilde b_1(z)^{[v+\ln |w|]_{q(z)}} F_2(z,\{v+\ln |w|\}_{q(z)})
1_{(-\infty,0)}(w)
$$
$$
+ F_3(z) [v+\ln |w|]_{q(z)} 1_{\{\widetilde b_1(z) = 1\}}1_{(-\infty,0)}(w),
$$
a.e.\ $\sigma(dz)dvdu$, for some functions $F_2$ and $F_3$. As in the case $\kappa \neq
0$ above, we can conclude that $X_\alpha$ can be represented by (\ref{e:representation})
where $\ln|t+u|/|u|$ in the last integrand term of (\ref{e:representation}) is replaced
by
$$
[v +  \ln |t+u| ]_{q(z)} - [v + \ln |u| ]_{q(z)}.
$$
Observe that, by using (\ref{e:int-frac-parts}), this difference can be expressed as
\begin{equation}\label{e:rewriting-last-term}
\frac{1}{q(z)} \Big( \{v +  \ln |t+u| \}_{q(z)} - \{v + \ln |u| \}_{q(z)} -
\ln\frac{|t+u|}{|u|}\Big).
\end{equation}
By including the first two terms of (\ref{e:rewriting-last-term}) into the first four
terms of (\ref{e:representation}), we can deduce that $X_\alpha$ can indeed be
represented by (\ref{e:representation}). \ \ $\Box$

\bigskip

{\sc Proof of Theorem \ref{t:identification-pfsm-canonical}:} Suppose that $X_\alpha$ is
represented by the sum of two independent processes (\ref{e:representation}) and
(\ref{e:mlfsm2}). To show that $X_\alpha$ is a PFSM, it is enough to prove that
(\ref{e:representation}) and (\ref{e:mlfsm2}) are PFSMs. The process
(\ref{e:representation}) has the representation (\ref{e:mma}) with
(\ref{e:cpfsm-X})--(\ref{e:cpfsm-G}). It is easy to verify that
$$
G(z,v,c(z)u) = c(z)^{\kappa} b_1(z) G(z,v,u) + F_3(z)q(z) 1_{\{b_1(z)=1\}}1_{\{\kappa
=0\}},
$$
where $c(z) = e^{q(z)}$. Hence, $C_P = Z\times [0,q(\cdot))$ where $C_P$ is the PFSM set
defined by (\ref{e:pfsm-set}). Definition \ref{d:pfsm-cfsm-kernel} yields that
(\ref{e:representation}) is a PFSM. One can show that the process (\ref{e:mlfsm2}) is a
PFSM in a similar way.

Suppose now that $X_\alpha$ is a PFSM. By Definition \ref{d:pfsm-cfsm-flow}, a minimal
representation of the process $X_\alpha$ is generated by a periodic flow. Suppose that
the process $X_\alpha$ has the minimal representation
$$
X_\alpha(t) \stackrel{d}{=} \int_{\widetilde X} \int_\bbR \widetilde G_t(\widetilde x,u)
\widetilde M_\alpha(d\widetilde x,du),
$$
where $(\widetilde X,\widetilde {\cal X},\widetilde \mu)$ is a standard Lebesgue space,
$\widetilde G_t(\widetilde x,u) = \widetilde G(\widetilde x,t+u) - \widetilde
G(\widetilde x,u)$, $\widetilde x\in \widetilde X,u\in\bbR$, $\{\widetilde
G_t\}_{t\in\bbR}\subset L^\alpha(\widetilde X\times \bbR,\widetilde\mu(dx)du)$, and
$\widetilde M_\alpha(d\widetilde x,du)$ has the control measure $\widetilde
\mu(d\widetilde x)du$, and that it is generated by a periodic flow $\{\widetilde
\psi_c\}_{c>0}$ on $\widetilde X$. Since the flow is periodic, we have $\widetilde X =
\widetilde P$, where $\widetilde P$ is the set (\ref{e:P}) of periodic points of the flow
$\{\widetilde \psi_c\}_{c>0}$. Partitioning $\widetilde P$ into the set $\widetilde L$ of
the cyclic points of the flow in (\ref{e:L}) and the set $\widetilde F$ of the fixed
points of the flow in (\ref{e:F}), we get
$$
X_\alpha(t) \stackrel{d}{=} \int_{\widetilde L} \int_\bbR \widetilde G_t(\widetilde x,u)
\widetilde M_\alpha(d\widetilde x,du) + \int_{\widetilde F} \int_\bbR \widetilde
G_t(\widetilde x,u) \widetilde M_\alpha(d\widetilde x,du) =: X_\alpha^L(t) +
X_\alpha^F(t).
$$
The processes $X_\alpha^L$ and $X_\alpha^F$ are independent since the sets $\widetilde L$
and $\widetilde F$ are disjoint. The process $X_\alpha^L$ is generated by a cyclic flow,
the flow $\{\widetilde \psi_c\}_{c>0}$ restricted to the set $\widetilde L$, and has a
representation (\ref{e:representation}) by Proposition \ref{p:representation-cyclic-flow}
below. The process $X_\alpha^F$ is generated by an identity flow satisfying $\widetilde
\psi_c(\widetilde x) = \widetilde x$ for all $\widetilde x\in \widetilde F$, and has a
representation (\ref{e:mlfsm}) by Theorem 5.1 in Pipiras and Taqqu
\cite{pipiras:taqqu:2002s}.
 \ \ $\Box$

\small


\noindent Vladas Pipiras \hfill  Murad S.\ Taqqu

\noindent  Department of Statistics and Operations Research    \hfill Department of
Mathematics and Statistics

\noindent  University of North Carolina at Chapel Hill \hfill Boston University

\noindent  CB\#3260, New West \hfill 111 Cummington St.

\noindent  Chapel Hill, NC 27599, USA \hfill Boston, MA 02215, USA

\noindent {\it pipiras@email.unc.edu} \hfill {\it murad@math.bu.edu}

\end{document}